\documentclass[11pt,reqno]{amsart}
\usepackage{mathrsfs}
\usepackage{amssymb}
\usepackage{amsthm}
\usepackage{amsmath,amsfonts,amssymb,esint}
\usepackage[dvips]{graphicx, color}

\newcommand{\nextpar}[1]{\vspace{10pt}\noindent\textit{#1}.\,}
\newcounter{counterpar}
\newcounter{ctr}
\newtheorem{definition}[ctr]{Definition}
\newtheorem{theorem}[ctr]{Theorem}
\newtheorem{proposition}[ctr]{Proposition}
\newtheorem{lemma}[ctr]{Lemma}

\newtheorem{corollary}[ctr]{Corollary}

\newcommand{\startproof}{\noindent{\bf Proof}.\,}
\newcommand{\stopproof}{\hfill$\blacksquare$\\}

\newcommand{\defn}[1]{\textbf{#1}}
\newcommand{\R}{{\mathbb R}}

\title{
h-principles for the incompressible Euler equations
}
\author{A.~Choffrut}
\date{\today}
\begin{document}
\maketitle

\begin{abstract}
In \cite{DS-Hoelder},
De Lellis and Sz\'ekelyhidi construct H\"older continuous,
dissipative (weak) solutions to the incompressible Euler equations 
in the torus $\mathbb T^3$.
The construction consists in adding fast oscillations to the trivial solution.
We extend this result by establishing optimal h-principles in two and three space dimensions.
Specifically, we identify all subsolutions (defined in a suitable sense)
which can be approximated in the $H^{-1}$-norm by exact solutions.
Furthermore, we prove that the flows thus constructed on $\mathbb T^3$ 
are genuinely three-dimensional and are not trivially obtained from solutions on $\mathbb T^2$.
\end{abstract}

%
%


\section{Introduction}
\subsection{Incompressible Euler equations and h-principle}
\setcounter{counterpar}{0}

We consider the \defn{(incompressible) Euler equations}
\begin{equation}\partial_tv+{\rm div}\,(v\otimes v)+\nabla p=0,\qquad{\rm div}\,v=0
\label{eq:Euler}
\end{equation}
on the torus $\mathbb T^d$, $d=2$ or $3$.
Here, $v$ is the velocity vector field and the pressure $p$ 
enforces the divergence-free condition.
If $(v,p)$ is a {\it classical} solution to \eqref{eq:Euler}, 
scalar multiplication with $v$ and the chain rule give
$\partial_t \frac{|v|^2}{2} + {\rm div}_x \left(\left(\frac{|v|^2}{2} + p\right) v \right) = 0$.
Integrating in space shows that classical solutions to the incompressible Euler equations 
conserve the total kinetic energy:
\begin{equation}
\frac{d}{dt}  \int_{\mathbb T^d} |v|^2 (x,t)\, dx = 0\,
\notag
\end{equation}

\nextpar{Anomalous dissipation}
The existence of \textit{weak} solutions violating 
the conservation of kinetic energy was first suggested in \cite{Onsager} by Onsager, 
where indeed he conjectured the existence of 
H\"older continuous solutions in $3$ space dimensions with any exponent smaller than $\frac{1}{3}$.
Onsager also asserted that such solutions do not exist if we impose the H\"older 
continuity with exponent larger than $\frac{1}{3}$ and this part of his conjecture was proved in \cite{Eyink} and
\cite{ConstantinETiti}. The considerations of Onsager are motivated by the Kolmogorov theory of isotropic
$3$-dimensional turbulence, where the phenomenon of anomalous dissipation in the Navier-Stokes equations
is postulated. This assumption seems to be widely confirmed experimentally, whereas no such
phenomenon is observed in $2$ dimensions. Indeed, for $d=2$ the conservation law for the enstrophy does
prevent it for solutions which start from sufficiently smooth initial data. 
However, the considerations put forward by Onsager which pertain to the mathematical
structure of the equations do not depend on the dimension and this independence appears clearly also in
the proof of \cite{ConstantinETiti}, which works for any $d\geq 2$.

\medskip

The first proof of the existence of a weak solution violating the energy conservation 
was given in 
the groundbreaking work of Scheffer \cite{Scheffer93}, 
which showed the existence of a compactly supported nontrivial
weak solution in $\R^2\times \R$. 
A different construction of the existence of a compactly supported nontrivial
weak solution in $\mathbb T^2\times \R$ was then given by Shnirelman in \cite{Shnirelman1}. 
In both cases the solutions are only square summable as a function of both space and time variables.
The first proof of the existence of a solution for which the total
kinetic energy is a monotone decreasing function has been
given by Shnirelman in \cite{Shnirelmandecrease}. 
Shnirelman's example is in the energy space $L^\infty ([0, \infty), L^2 (\R^3))$. 

In \cite{DS1,DS2} these existence results were extended to solutions 
with bounded velocity and pressure and in any space dimensions. 
The same methods were also used to give quite severe counterexamples to the uniqueness
of admissible solutions, both for incompressible and compressible Euler. 
Further developments in fluid dynamics inspired by these works appeared 
subsequently in \cite{Chiodaroli, CFG, Shvydkoy, Szekelyhidi, SzWie, Wiedemann} 
and are surveyed in the note \cite{DSsurvey}. 
In \cite{DS3, DS-Hoelder}, De Lellis and Sz\'ekelyhidi devised a new iteration 
scheme , which produces continuous and even H\"older continuous solutions
on $\mathbb T^3$.
Furthermore, one may prescribe the \defn{total kinetic energy} profile
\begin{equation}\fint_{\mathbb T^d} |v(x,t)|^2\,dx=e(t)\notag
\end{equation}
where $d=3$, and $\fint_{\mathbb T^d}=\frac1{(2\pi)^d}\int_{\mathbb T^d}$.
(For notational convenience we omit the usual factor $1/2$ and average over the domain.)
\medskip

Solutions of class $C^1$ are therefore ``rigid'' compared to less regular solutions.
In fact, the paper \cite{DS1} introduced a new point of view in the subject,
highlighting connections to other counterintuitive solutions of 
(mainly geometric) systems of partial differential equations: 
in geometry these solutions are, according to Gromov, instances of the $h$-principle,
the prime example of which is Nash's theorem on $C^1$ isometric embeddings \cite{Nash54}.
(See in \cite{CDSz} an earlier discussion on the striking similarities between
Onsager's conjecture and 
the rigidity and flexibility properties of the isometric problem.)
We recall that an embedding $u_0\colon M^n\rightarrow \mathbb R^N$, $N>n$ 
is said to be (strictly) short
if $\partial_iu_0\cdot \partial_ju_0< g_{ij}$ where $g$ is a prescribed Riemannian metric.
Nash (and Kuiper) proved that any strictly short embedding can be uniformly approximated 
by an isometric embedding $u\in C^1(M;\mathbb R^N)$, $\partial_iu\cdot\partial_ju=g_{ij}$,
in the sense that $\|u_0-u\|_{C^0}$ can be made arbitrarily small.
For the isometric problem, the h-principle is the statement 
that $u_0$ can be deformed into $u$ {\it via} a \textit{homotopy} (hence the name).
In the sequel we will leave this aspect of the h-principle aside
and view the h-principle as a density statement.

The main idea in \cite{Nash54} is to add fast oscillations in order to increase 
the metric induced by a short embedding $u_0$ 
and thereby reducing the defect $g_{ij}-\partial_iu_0\cdot \partial_ju_0$.
Thus, $u_0$ is taken closer to the boundary of the set short embeddings,
precisely made up of isometric embeddings.
Nash's idea has been further developped by Kuiper, Gromov, and others,
and falls nowadays under the name of \defn{convex integration}, 
see \cite{DSsurvey, Eliashberg, Gromov, SPRING:book}. 

If convex integration alone produces $C^1$ isometric embeddings,
refinements can achieve $C^{1,\alpha}$ regularity for certain $\alpha$ depending on $n$ and $N$,
see \cite{Borisov65, Borisov2004, CDSz} for precise statements and references therein.
For the Euler equations, the natural space for convex integration is $C^0$.
The method used in \cite{DS1} producing solutions in $L^\infty$ was a weak form of convex integration.
The iteration scheme of \cite{DS3} is closer to the approach of \cite{Nash54},
see the introduction of \cite{DS3} for a thorough discussion.
Finally, \cite{DS-Hoelder}, with the improved regularity for Euler, parallels \cite{CDSz} for the isometric problem.
\medskip

In this article we establish h-principles for the Euler equations in $2$ and $3$ space dimensions,
using the convex integration procedure developped in \cite{DS3} and sophisticated in \cite{DS-Hoelder}.
We shall first motivate our definition of subsolutions to the Euler equations,
analogous to the short embeddings of Nash for the isometric problem,
and the notion of the h-principle in use here.
It is generally accepted that the onset of turbulence in incompressible fluids is due to
the appearance of high-frequency oscillations in the velocity field
\cite{DSsurvey, FrischBook, MAJDA:nonlinear-analysis-applied-math}.
For example, if $(v_\nu, p_\nu,f_\nu)$ is a sequence of approximate solutions, 
$$\partial_tv_\nu+{\rm div}\,(v_\nu\otimes v_\nu)+\nabla p_\nu=f_\nu,\qquad {\rm div}\,v_\nu=0,$$
with uniformly bounded kinetic energy, 
and converges weakly to $(v,\pi,0)$,
then in general
\begin{equation}\partial_tv+{\rm div}\,(v\otimes v+R)+\nabla \pi=0,\qquad {\rm div}\,v=0
\label{eq:system-in-(v,R,pi)}
\end{equation}
in the weak sense, where $R$ is a symmetric, positive semi-definite matrix,
called the \defn{Reynolds stress tensor}.
It appears because the operation of taking weak limits does not commute
with the nonlinear operator $\otimes$.
A strategy to construct an exact solution to the Euler equations (\ref{eq:Euler})
is then to reintroduce the oscillations so as to eliminate $R$ on average.
A crucial point in the construction of \cite{DS3} is therefore the ability
to generate the tensor $R\geq 0$ with a fast oscillating perturbation $W$
(see Lemma~\ref{lemma:geometric-lemma} and Section~\ref{section:definitions-of-the-maps} for details):
we seek a velocity field $W$ solving the stationary
Euler equations and satisfying
$$\fint_{\mathbb T^d}W\otimes W\,d\xi=R.$$
In three dimensions, this is done using Beltrami flows.
However, these flows seem to be insufficient to capture all possible oscillatory behaviors
in the Euler equations,
see Proposition~\ref{proposition:geometry-of-sym-pos-def-matrices},
where it is also shown that this problem does not exist in two dimensions.
\medskip

\noindent{\bf Remark}\quad
Beltrami flows are defined as those three-dimensional flows of the form
${\rm curl}\,v(x)=\lambda(x)v(x)$ for some scalar function $\lambda(x)$.
These are stationary flows, see \cite{MajdaBook}.
There is a connection with two-dimensional stationary flows,
see Proposition~2.11 in \cite{MajdaBook}.
For these flows however, the function $\lambda$ is in general not constant.
In the construction of \cite{DS3}, on the other hand,
the function $\lambda(x)=\lambda$ is constant.
\stopproof

With these general considerations being done, we now turn to precise definitions.
It will be more convenient to work with an alternative form of (\ref{eq:system-in-(v,R,pi)}).
Letting $\mathring R$ be (minus) the trace-free part of $R$, 
$R=\frac{{\rm tr}\,R}d\,{\rm Id}-\mathring R$,
then (\ref{eq:system-in-(v,R,pi)}) becomes
\begin{equation}
\partial_tv+{\rm div}\,(v\otimes v)+\nabla p={\rm div}\,\mathring R,\qquad {\rm div}\,v=0
\label{eq:Euler-Reynolds}
\end{equation}
where $p=\pi+\frac{{\rm tr}\,R}d$.
We shall refer to (\ref{eq:Euler-Reynolds}) as the \defn{Euler-Reynolds system}.
It is equivalent to (\ref{eq:system-in-(v,R,pi)}) provided one fixes ${\rm tr}\,R$.
(Indeed, if $(v,R,\pi)$ solves (\ref{eq:system-in-(v,R,pi)}),
then so does $(v,R+f{\rm Id},\pi-f)$ for any function $f$.)

We emphasize that the notion of short embedding for the isometric problem
is relative to a prescribed metric $g$.
In the context of ideal hydrodynamics,
a natural quantity to prescribe is the kinetic energy $e(t)$.
We shall say that $(v,\mathring R,p)$ is a \defn{strict subsolution} to the Euler equations 
(relative to the kinetic energy $e(t)$)
if $(v,\mathring R,p)$ solves the Euler-Reynolds system (in the classical sense),
where $\mathring R$ is trace-free, and if 
\begin{equation}\frac{e-\fint_{\mathbb T^d}|v(x',t)|^2\,dx'}d{\rm Id}-\mathring R(x,t)>0,
\qquad x\in\mathbb T^d,\quad t\in[0,T].
\label{ineq:subsolution-in-terms-of-(v,Ro,p)}
\end{equation}
This amounts to fixing ${\rm tr}\,R=e(t)-\fint_{\mathbb T^d}|v(x,t)|^2\,dx$ 
and imposing that $R>0$ in the $(v,R,\pi)$ formulation.
In particular we have 
\begin{equation}e>\fint_{\mathbb T^d}|v|^2.\notag
\end{equation}
As for the isometric problem, the boundary of the set of subsolutions consists of exact solutions
to the Euler equations with prescribed kinetic energy.

The main focus \cite{DS3, DS-Hoelder} is the construction of
{\it some} solutions with a certain amount of regularity,
and thus used the particular (trivial) subsolution $(0,0,0)$.
Their Geometric Lemma (Lemma~3.2 in \cite{DS3}) was sufficient for this purpose.
Here, we prove an optimal Geometric Lemma, see Lemma~\ref{lemma:geometric-lemma},
and identify the largest class of subsolutions for which the convex integration scheme of \cite{DS3}
produces an exact solution to (\ref{eq:Euler}).
A subsolution $(v,\mathring R,p)$ is \defn{strong} if it 
satisfies the condition, stronger than (\ref{ineq:subsolution-in-terms-of-(v,Ro,p)}),
that
\begin{equation}
\frac{e(t)-\fint_{\mathbb T^d}|v_0(x',t)|^2\,dx'}{d(d-1)}{\rm Id}+\mathring R_0(x,t)>0,
\qquad x\in\mathbb T^d,\quad t\in[0,T].
\label{ineq:strong-subsolution}
\end{equation}
(Equivalently,
$\frac{e(t)-\fint_{\mathbb T^d}|v_0(x',t)|^2\,dx'}d{\rm Id}-\mathring R_0(x,t)\in \mathcal M_d$,
see Section~\ref{section:notation} for definitions.)
\medskip

We say that the \defn{h-principle} holds for (\ref{eq:Euler})
if, 
given $\sigma>0$ and a strict subsolution $(v_0,\mathring R_0,p_0)$ relative to $e(t)$,
there exists an exact solution $(v,p)$ with $\fint_{\mathbb T^d}|v(x,t)|^2\,dx=e(t)$
and such that $\|v-v_0\|_{H^{-1}(\mathbb T^d)}<\sigma$.
The h-principle holds \defn{for strong subsolutions}
if $(v_0,\mathring R_0,p_0)$ is required to be strong.
\medskip

\noindent{\bf Remark}\quad
Another possible notion of subsolutions is to fix a function $\overline e=\overline e(x,t)$
and impose the pointwise condition that
$\frac{\overline e(x,t)-|v(x,t)|^2}d\,{\rm Id}-\mathring R(x,t)> 0$.
This is the notion used in \cite{DS1} in the context of $L^\infty$-solutions.
The two notions are different, and one does not imply the other.
The pointwise notion seems ill suited for the construction in use here.
Indeed, the pointwise control on the velocity field along the iteration seems insufficient.
\medskip

\nextpar{Genuinely 3D flows}
There is a trivial way to produce flows on $\mathbb T^3$ from flows
on $\mathbb T^2$.
As a consequence of precise estimates of our main result,
Theorem~\ref{theorem:weak-h-principle},
we show that the flows obtained for $d=3$ 
are genuinely three-dimensional and do not coincide with those obtained for $d=2$,
see Corollary~\ref{corollary:3D-flows-are-not-2D}.
In order to formulate our statement precisely,
consider a solution $(v,p)$ to the Euler equations (\ref{eq:Euler})
on $\mathbb R^3\times [0,T]$
and denote with the same letters
the corresponding solution on $\mathbb R^3\times [0,T]$
with the obvious periodicity in space.
We say that a solution is \defn{not genuinely three-dimensional} if,
after suitably changing  coordinates in space,
it takes the form
\begin{equation}v(x,t)=(v_1(x_1,x_2,t),v_2(x_1,x_2,t),v_3)
\label{eq:3D-flow-which-is-2D}
\end{equation}
where $v_3$ is a constant.
Otherwise it is \defn{genuinely three-dimensional}.
\bigskip

\noindent{\bf Acknowledgments}\quad
The author was supported by ERC Grant agreement No.~277993.
The author is extremely grateful to Camillo De Lellis and L{\'a}szl{\'o} Sz{\'e}kelyhidi, Jr.
for introducing him to this subject and for fruitful discussions.

\subsection{Statement of results}
In \cite{DS3} and \cite{DS-Hoelder}, solutions to the Euler equations with prescribed
kinetic energy were constructed using convex integration starting from the trivial 
subsolution $(v_0,\mathring R_0,p_0)=(0,0,0)$ on $\mathbb T^3$.
Since the building blocks are a certain class of Beltrami flows, which are inherently three-dimensional,
it is not immediately clear whether the method should work in other space dimensions.
In our main result, Theorem~\ref{theorem:weak-h-principle}, 
we establish the largest set of subsolutions for which the h-principle holds,
in dimensions two and three.
It is based on a refined Geometric Lemma,
see Proposition~\ref{proposition:geometry-of-sym-pos-def-matrices}
and Lemma~\ref{lemma:geometric-lemma}.
\begin{theorem}[h-principle]\label{theorem:weak-h-principle}
Assume $d=2$ or $3$.
Let $e(t)$, $t\in[0,T]$, be smooth, positive.
Let $(v_0,\mathring R_0,p_0)$ be a strong subsolution. 
Let $0<\theta<\frac1{10}$ and $\sigma>0$.
Then:
\begin{enumerate}
\item there exists a vector field $v\in C^0(\mathbb T^d\times [0,T])$
and a function $p\in C^0(\mathbb T^d\times [0,T])$
which solve the Euler equations (\ref{eq:Euler}) (in the weak sense) and satisfy
$$|v(x,t)-v(x',t)|\leq C|x-x'|^\theta\qquad x,x'\in\mathbb T^d,\qquad t\in[0,T]$$
and
\begin{equation}\sup_{t\in[0,T]}\|v(\cdot,t)-v_0(\cdot,t)\|_{H^{-1}(\mathbb T^d)}<\sigma;
\label{ineq:|v|H-1<epsilon}
\end{equation}
\item the solution can be constructed so that, for all $t\in[0,T]$,
\begin{equation}
\left|\fint_{\mathbb T^d}\left(v(x,t)\otimes v(x,t)-v_0(x,t)\otimes v_0(x,t)+\mathring R_0(x,t)\right)\,dx
-\frac{e(t)-\fint_{\mathbb T^2}|v_0(x,t)|^2\,dx}d{\rm Id}\right|<\sigma.
\label{ineq:int(vxv)}
\end{equation}
\end{enumerate}
\end{theorem}

\noindent{\bf Remark}\quad
As in \cite{DS-Hoelder}, the proof of Theorem~\ref{theorem:weak-h-principle}
yields further regularity on both $v$ and $p$.
Namely, they are H\"older continuous in both $x$ and $t$ with
$$|v(x,t)-v(x',t')|\leq C\left(|x-x'|^\theta+|t-t'|^\theta\right)$$
and
$$|p(x,t)-p(x',t')|\leq C\left(|x-x'|^{2\theta}+|t-t'|^{2\theta}\right).$$
\stopproof

\begin{corollary}[]\label{corollary:strong-h-principle-for-2D-Euler}
If $d=2$, then the h-principle holds for strict subsolutions.
\end{corollary}


\begin{corollary}[Genuine 3D flows]\label{corollary:3D-flows-are-not-2D}
Assume $d=3$ in Theorem~\ref{theorem:weak-h-principle}.
Then the flows are genuinely three-dimensional provided $\sigma$ is chosen sufficiently small.
\end{corollary}

\section{Proof of Theorem~\ref{theorem:weak-h-principle}, part~1}
\subsection{Notation}\label{section:notation}
\setcounter{counterpar}{0}
\smallskip

\nextpar{Spaces of symmetric matrices}
All matrices in this article will be symmetric, and thus the qualifier ``symmetric'' will often be omitted.
We shall denote by 
$$\mathcal S^{d\times d}=\left\{M\in\mathbb R^{d\times d}~\colon~M^\top=M\right\}$$
the set of symmetric $d\times d$ matrices,
by
$$\mathcal S^{d\times d}_{++}=\left\{M\in\mathcal S^{d\times d}~\colon~M>0\right\}$$
the open convex cone of (symmetric) positive definite matrices,
and by
$$\mathcal  S^{d\times d}_0=\left\{M\in\mathcal S^{d\times d}~\colon~{\rm tr}\,M=0\right\}$$
the closed linear subset of $\mathcal S^{d\times d}$ of trace-free matrices.
We also introduce
$$\mathbb M_d:=\left\{{\rm Id}-b\otimes b~|~b\in\mathbb S^{d-1}\right\}$$
where $\mathbb S^{d-1}=\left\{b\in \mathbb R^d~|~|b|=1\right\}$ denotes the $(d-1)$-dimensional sphere.
Of interest will be the open subset 
$$\mathcal M_d:={\rm int}\,\mathbb M_d^{\rm conic}\subset \mathcal S^{d\times d}_{++}
$$
where $\mathbb M_d^{\rm conic}$ denotes the \defn{conic hull} of $\mathbb M_d$,
that is, the set of all matrices of the form
$$\sum_{i=1}^m\alpha_i\left({\rm Id}-b_i\otimes b_i\right),\qquad{\rm where}\quad  \alpha_i>0\quad{\rm and}\quad |b_i|=1.$$

The norm on these spaces will be the operator norm.

\nextpar{$B_{r_o}({\rm Id})$ and the parameter $\overline r_0$}
For $r_0>0$, $B_{r_0}({\rm Id})$ will always denote the open ball in $\mathcal S^{d\times d}$.
By Proposition~\ref{proposition:geometry-of-sym-pos-def-matrices},
we fix $\overline r_0>0$ sufficiently small so that
\begin{equation}\overline{B_{2\overline r_0}({\rm Id})}\subset \mathcal M_d.
\label{eq:fix-r0-bar}
\end{equation}

\nextpar{H\"older norms}
For a time-independent function $f=f(x)$, 
the sup-norm is denoted $\|f\|_0=\sup_{\mathbb T^d}|f|$, and 
the H\"older seminorms are given by
\begin{eqnarray}
\left[f\right]_m:=\max_{|\gamma|=m}\|D^\gamma f\|_0,\qquad
\left[f\right]_{m+\alpha}:=\max_{|\gamma|=m}\sup_{x\neq y}\frac{|D^\gamma f(x)-D^\gamma f(y)|}{|x-y|^\alpha}
\notag
\end{eqnarray}
and the H\"older norms are given by
\begin{eqnarray}
\|f\|_m:=\sum_{j=0}^m\left[f\right]_j,\qquad
\|f\|_{m+\alpha}:=\|f\|_m+\left[f\right]_{m+\alpha}.
\notag
\end{eqnarray}

For a time-dependent function $f=f(x,t)$, and $r\geq 0$,
$\|f\|_r$ will denote the ``H\"older norm in space'', that is
\[
\|f\|_r = \sup_{t\in [0,T]} \|f (\cdot, t)\|_r\, .
\]
while the H\"older norms in space and time will be denoted by $\|\cdot\|_{C^r}$.

\nextpar{Constants}
We will follow \cite{DS-Hoelder} for the convention pertaining to the constants
involved in the estimates of Section~\ref{section:proof-of-Proposition(iteration)} and the Appendix.
\begin{itemize}
\item $C$: will denote universal constants.
\item $C_h$: will denote constants in estimates concerning standard functional inequalities in H\"older spaces $C^r$.
These constants depend only on the specific norm used and therefore only on the parameter $r\geq 0$.
\item $C_e$: throughout the rest of the paper the prescribed energy will be assumed to be a fixed
smooth function bounded below by a positive function.
Several estimates depend on these bounds and the relate constants will be denoted $C_e$.
\item $C_v$: in addition to the dependence on $e$, there will be estimates which depend also on $\|v\|_0$.
See the constant $A$ in Proposition~\ref{proposition:iteration}.
\item $C_s$, $C_{e,s}$, $C_{v,s}$: will denote constants which are typically involved
in Schauder estimates for $C^{m+\alpha}$ norms of elliptic operators,
when $m\in\mathbb N$ and $0<\alpha<1$.
These constants not only depend on the specific norm used, but they also degenerate as $\alpha\downarrow 0$ and $\alpha\uparrow 1$.
The ones denoted by $C_{e,s}$ and $C_{v,s}$ depend also, respectively, on $e$ and $e$ and $\|v\|_0$.
\end{itemize}
We emphasize that constants never depend on the parameters $\mu,\ell, \delta,\lambda$ and $D$, 
although they may depend on $\varepsilon$
(see Section~\ref{section:proof-of-Proposition(iteration)} and Appendix
for definitions of these parameters).

\subsection{The iterative scheme: Proposition~\ref{proposition:iteration}}
In order to motivate the main Proposition of this Section,
we briefly sketch the strategy to construct exact solutions to the Euler equations (\ref{eq:Euler}).
Given a strict subsolution $(v,\mathring R,p)$,  {\it i.e.}
$$\partial_tv+{\rm div}\,(v\otimes v)+\nabla p={\rm div}\,\mathring R,\qquad {\rm div}\,v=0$$
and
$$\frac{e(t)-\fint_{\mathbb T^d}|v(x,t)|^2\,dx}d\,{\rm Id}-\mathring R>0,$$
we construct a triple $(v_1,\mathring R_1,p_1)$ which is closer to being a solution,
in the sense that the energy gap $e(t)-\fint_{\mathbb T^d}|v_1(x,t)|^2\,dx$ 
and the trace-free tensor $\mathring R_1$ are both smaller, while
$\frac{e(t)-\fint_{\mathbb T^d}|v_1(x,t)|^2\,dx}d\,{\rm Id}-\mathring R_1$ remains positive definite.
An iteration is needed since $e(t)-\fint_{\mathbb T^d}|v_1(x,t)|^2\,dx$ 
and $\mathring R_1$ cannot be made to vanish exactly.
Yet the iteration converges because this can be done with arbitrary accuracy.
\begin{proposition}\label{proposition:iteration}
Suppose $d=2$ or $3$, and fix $\overline r_0>0$ as in (\ref{eq:fix-r0-bar}).
Let $K\subset \mathcal M_d$ be compact and contain $B_{2\overline r_0}({\rm Id})$,
and let $\mathcal N$ be an open neighborhood of $K$ such that $\overline{\mathcal N}\subset \mathcal M_d$.
Fix $e_0\geq 0$, $\Delta_0>0$,
and $\tilde r_0=\tilde r_0(d,K,e_0,\Delta_0)$ as in Lemma~\ref{lemma:wo-is-well-defined}.

Fix now $r_0\leq \min\{\overline r_0,\tilde r_0\}$ and set
\begin{equation}\eta = \frac{\min\Delta_0}{4d}r_0.
\label{eq:eta}
\end{equation}
Then, there exists $M=M(e_0,\Delta_0)$ with the following properties.
\bigskip

Let $\varepsilon>0$ and $0<\zeta\leq \frac12$.
Suppose $0<\delta\leq 1$ and $(v,\mathring R,p)$ satisfy
$$\partial_tv+{\rm div}\,(v\otimes v)+\nabla p={\rm div}\,\mathring R,\qquad {\rm div}\,v=0,$$
\begin{equation}
\left|e_0(t)+\Delta_0(t)(1-\delta)-\fint_{\mathbb T^d}|v(x,t)|^2\,dx\right|\leq \zeta/2\delta \Delta_0(t),
\label{ineq:assumed-energy-gap}
\end{equation}
and, posing $\overline\delta:=\zeta\delta^\frac32$, that
\begin{equation}
{\rm Id}-\frac d{e_0(t)+(1-\overline\delta)\Delta_0(t)-\fint_{\mathbb T^d}|v(x,t)|^2\,dx}\mathring R(x,t)\in K,\qquad x\in\mathbb T^d,\quad t\in\mathbb S^1.
\label{eq:rescaled-R-is-in-K}
\end{equation}
Set
$D:=\max\left\{1,\|v\|_1,\|\mathring R\|_1\right\}$.

Then there exists $(v_1,\mathring R_1,p_1)$ satisfying
$$\partial_tv_1+{\rm div}\,(v_1\otimes v_1)+\nabla p_1={\rm div}\,\mathring R_1,\qquad {\rm div}\,v_1=0$$
and such that
\begin{equation}
\left|e_0(t)+\Delta_0(t)(1-\overline\delta)-\fint_{\mathbb T^d}|v_1(x,t)|^2\,dx\right|
\leq \zeta/2\overline\delta \Delta_0(t),
\label{ineq:target-energy-gap}
\end{equation}
\begin{eqnarray}
\|\mathring R_1\|_0&\leq& \eta\overline\delta,\label{ineq:target-bound-on-|Ro1|0}\\
\|v_1-v\|_0&\leq& M\sqrt\delta,\label{ineq:target-bound-on-|v1-v|0}\\
\sup_t\|v_1(\cdot,t)-v(\cdot,t)\|_{H^{-1}(\mathbb T^d)}&\leq& r_0\delta^\frac14,\label{ineq:target-bound-on-|v1-v|H(-1)}\\
\|p_1-p\|_0&\leq& M^2\delta,\label{ineq:target-bound-on-|p1-p|0}\\
\max\left\{1,\|v_1\|_1,\|\mathring R_1\|_1\right\}&\leq& A\delta^\frac32\left(\frac D{\overline\delta^2}\right)^{1+\varepsilon}
\label{ineq:target-bound-on-D}
\end{eqnarray}
where the constant $A$ depends on $d, e, \varepsilon>0$ and $\|v\|_0$,
see (\ref{eq:A}).
\end{proposition}
\bigskip

\noindent{\bf Remark}\quad 
The conclusions imply that
\begin{equation}
{\rm Id}-\frac d{e_0(t)+(1-\overline{\overline\delta})\Delta_0(t)-\fint_{\mathbb T^d}|v_1(x,t)|^2\,dx}\mathring R_1(x,t)
\in B_{r_0}({\rm Id})\subset K,\qquad x\in\mathbb T^d,\quad t\in\mathbb S^1
\label{eq:rescaled-R1-is-in-K}
\end{equation}
where $\overline{\overline\delta}=\zeta\overline\delta^\frac32$.
Indeed,
\begin{eqnarray}
e_0+\Delta_0(1-\overline{\overline\delta})-\fint_{\mathbb T^d}|v_1|^2
=\Delta_0(\overline\delta-\overline{\overline\delta})
+e_0+\Delta_0(1-\overline\delta)-\fint_{\mathbb T^d}|v_1|^2
\geq \frac{\Delta_0\overline\delta}4
\notag
\end{eqnarray}
so that
$\left\|\frac d{e_0+\Delta_0(1-\overline{\overline\delta})-\fint_{\mathbb T^d}|v_1(x,t)|^2\,dx}\mathring R_1\right\|
\leq r_0.$
Therefore, an iteration can be carried out by repeated use of Proposition~\ref{proposition:iteration}.
\stopproof

\subsection{Proof of Theorem~\ref{theorem:weak-h-principle}, part 1 }
\label{section:proof-of-part-1-of-theorems}
\setcounter{counterpar}{0}
Assume Proposition~\ref{proposition:iteration} is proved.

\nextpar{Periodicity in $t$}
In this paragraph we show that we may assume, without loss of generality,
that $v_0,\mathring R_0,p_0$, and $e$ are periodic in $t$.
(Although this is not necessary for the construction,
this feature will prove to be convenient as mollification in space {\it and} time is used in the estimates,
see Section~\ref{section:definitions-of-the-maps}).
Let's then start with a strong subsolution $(v_0, \mathring R_0, p_0)$
defined for $x\in\mathbb T^d$ and $t\in[0,T]$
relative to $e(t)$ which is a smooth positive function defined for $t\in[0,T]$.
It is standard that $v_0(x,t), p_0(x,t)$ can be extended to smooth functions
for $x\in\mathbb T^d$ and $t\in\mathbb R$ which vanish for $t\geq \frac32T$ and $t\leq -\frac T2$,
and such that ${\rm div}\,v_0=0$ and $\fint_{\mathbb T^d}v_0(x,t)\,dx=0$ for all $t\in \mathbb R$.
See for instance the proof of Corollary~1.3.7, p.~138, Part~II of \cite{Hamilton:Nash-Moser-IFT}.
We may then repeat $v_0$ and $p_0$ periodically in $t$ with period $2T$.
We define
$$\mathring R_0:=\mathcal R\,\left(\partial_tv_0+{\rm div}\,(v_0\otimes v_0)+\nabla p_0\right)$$
for $t\not\in[0,T]$,
see Definition~\ref{definition:the-operator-R} for the operator $\mathcal R$.
Since the argument of the right-hand side has average $0$ over $\mathbb T^d$, 
the triple $(v_0,\mathring R_0,p_0)$ solves the Euler-Reynolds system and is periodic in $t$,
see Lemma~\ref{lemma:R=divInverse}.
Finally, it is clear that $e(t)$ can be extended to a smooth, positive, periodic functions 
for $t\in \mathbb R$ with period $2T$ such that
$$\frac{e(t)-\fint_{\mathbb T^d}|v_0(x',t)|^2\,dx'}d\,{\rm Id}-\mathring R_0(x,t)\in\mathcal M_d,
\qquad x\in\mathbb T^d,\quad t\in\mathbb R.$$
Rescaling in $t$, we may assume that $v_0,\mathring R_0, p_0$, and $e$ 
have period $2\pi$.

\nextpar{Setting parameters}
Set
$$e_0(t):=\fint_{\mathbb T^d}|v_0(x,t)|^2\,dx,\qquad \Delta_0(t):=e(t)-e_0(t).$$
We may then choose $0<\zeta\leq \frac12$ such that
\begin{equation}
\frac{(1-\zeta)\Delta_0(t)}d{\rm Id}-\mathring R_0(x,t)\in\mathcal M_d,\qquad x\in\mathbb T^d,\quad t\in\mathbb R.
\label{ineq:condition-on-zeta-so-that-initial-R-can-be-generated}
\end{equation}
Set 
$$K:=\left\{{\rm Id}-\frac d{(1-\zeta)\Delta_0(t)}\mathring R_0(x,t)~\colon~ x\in\mathbb T^d, t\in[0,T]\right\}\cup \overline{B_{2\overline r_0}({\rm Id})},$$
where $\overline r_0$ is as in (\ref{eq:fix-r0-bar}).
Fix an open neighborhood $\mathcal N$ of $K$ such that $\overline{\mathcal N}\subset \mathcal M_d$.
Fix $\tilde r_0$ as in Lemma~\ref{lemma:wo-is-well-defined}
and let $r_0\leq \min\{\overline r_0,\tilde r_0\}$ to be specified later,
see Section~\ref{section:fixing-the-parameters}.

Define inductively
\begin{equation}\delta_{n+1}=\zeta\delta_n^\frac32,\qquad \delta_0=1.
\notag
\end{equation}

Fix $\varepsilon>0$ and $\sigma>0$.

\nextpar{The iterates}
Use Proposition~\ref{proposition:iteration} inductively to construct a sequence $(v_n,\mathring R_n,p_n)$
with $\delta=\delta_n$, $\overline \delta=\delta_{n+1}$.
Since $(v_0,\mathring R_0,p_0)$ clearly satisfies 
the assumptions of Proposition~\ref{proposition:iteration},
the remark following Proposition~\ref{proposition:iteration} shows that
${\rm Id}-\frac d{e_0+\Delta_0(1-\delta_{n+1})-\fint_{\mathbb T^d}|v_n|^2}\mathring R_n\in K$ for each $n$.
The sequence satisfies
\begin{equation}
\left|e_0(t)+\Delta_0(t)(1-\delta_n)-\fint_{\mathbb T^d}|v_n(x,t)|^2\,dx\right|\leq \zeta/2\delta_n \Delta_0(t),
\label{ineq:n-th-step-energy-gap}
\end{equation}
for $n\geq 0$, and for $n\geq 1$
\begin{eqnarray}
\|\mathring R_n\|_0&\leq& \eta\delta_n,\label{ineq:n-th-step-bound-on-Ro}\\
\|v_n-v_{n-1}\|_0&\leq& M\sqrt\delta_{n-1},\label{ineq:n-th-step-bound-on-(v(n)-v(n-1))}\\
\sup_t\|v_n(\cdot,t)-v_{n-1}(\cdot,t)\|_{H^{-1}(\mathbb T^d)}&\leq& r_0\delta_{n-1}^\frac14
\label{ineq:n-th-step-bound-on-H(-1)-norm}\\
\|p_n-p_{n-1}\|_0&\leq& M^2\delta_{n-1},\label{ineq:n-th-step-bound-on-(p(n)-p(n-1))}
\end{eqnarray}
and
\begin{equation}
D_{n+1}:=\max\left\{\|v_{n+1}\|_{C^1},\|\mathring R_{n+1}\|_{C^1}\right\}
\leq A\delta_n^\frac32\left(\frac{D_n}{\delta_{n+1}^2}\right)^{1+\varepsilon}
\label{ineq:n-th-step-bound-on-D(n+1)}
\end{equation}

\nextpar{Convergence of $\delta_n$ and $D_n$}
With $d_n=\ln(\zeta^2\delta_n)$ we have $d_{n+1}=\frac32d_n$ and so
\begin{equation}
\delta_n=\zeta^{-2}\zeta^{3(\frac32)^{n-1}}\qquad (n\geq 0).
\notag
\end{equation}

Next, define 
$x_n:=\delta_n^\gamma D_n$
where $\gamma>0$ will be chosen later.
Then, (\ref{ineq:n-th-step-bound-on-D(n+1)}) gives
\begin{eqnarray}
x_{n+1}&\leq& A\zeta^{-2(1+\varepsilon)+\gamma}\,\delta_n^{\gamma(\frac12-\varepsilon)-3(\frac12+\varepsilon)}\,x_n^{1+\varepsilon}.
\notag
\end{eqnarray}
There is no loss in assuming that $\varepsilon<\frac12$ (since we will take $\varepsilon\downarrow 0$).
Let
$$\gamma>3\frac{1+2\varepsilon}{1-2\varepsilon}$$
and observe that $0<\delta_n\leq 1$ so that
$x_{n+1}\leq A\zeta^{-2(1+\varepsilon)+\gamma}\,x_n^{1+\varepsilon}$.
Let $B:=\left(A\zeta^{-2(1+\varepsilon)+\gamma}\right)^{-\frac1\varepsilon}$, pose $z_n=\ln(Bx_n)$
One easily finds
$$D_{n+1}
=\zeta^{2\gamma}B^{-1}\zeta^{-3\gamma(\frac32)^{n+1}}(Bx_1)^{(1+\varepsilon)^n}.$$
Since $0<\varepsilon<\frac12$, $\gamma>0$, and $\ln\zeta<0$,
the term in $(\frac32)^n$ will dominate that in $(1+\varepsilon)^n$.
That is, for any $\gamma'>\gamma$, there exists $C'=C'(\zeta,\varepsilon,A,\gamma,\gamma')$ such that
$$D_{n+1}\leq C'\zeta^{-3\gamma'(\frac32)^n}.$$

\nextpar{Convergence in $C^0$ and (weak) solution to the Euler equations}
Since $\delta_n$ vanishes very fast,
and from (\ref{ineq:n-th-step-energy-gap}), (\ref{ineq:n-th-step-bound-on-Ro}), 
(\ref{ineq:n-th-step-bound-on-(v(n)-v(n-1))}), and (\ref{ineq:n-th-step-bound-on-(p(n)-p(n-1))}),
we conclude that $(v_n,p_n)$ converges uniformly to a (weak) solution $(v,p)$
to the Euler equations (\ref{eq:Euler})
with kinetic energy $e(t)=\fint_{\mathbb T^d}|v(x,t)|^2\,dx$.
In fact,
\begin{equation}
\|v_n-v_0\|_0\leq M\sum_{j=0}^\infty \delta_n^\frac12\leq CM
\label{ineq:a-priori-bound-on-|v(n)-v(0)|}
\end{equation}
where $C$ is some universal constant.
In turn, the constant $A$ in Proposition~\ref{proposition:iteration} can be taken 
to depend only on $\varepsilon$ and $e$.

\nextpar{Convergence in $C^\theta$}
We have
$$\|v_{n+1}-v_n\|_0\leq M\sqrt{\delta_n}\leq M\zeta^{-1}\zeta^{(\frac32)^n}$$
$$\|v_{n+1}-v_n\|_{C^1}\leq D_n+D_{n+1}\leq C'\zeta^{-3\gamma'(\frac32)^n}$$
and therefore by interpolation we find
$$\|v_{n+1}-v_n\|_{C^{\theta}}\leq \|v_{n+1}-v_n\|_{C^0}^{1-\theta}\|v_{n+1}-v_n\|_{C^1}^\theta
\leq (M\zeta^{-1})^{1-\theta}(C')^\theta\,\zeta^{\left((1-\theta)-3\gamma'\theta\right)(\frac32)^n}.$$
The critical value for $\theta$ for which the right-hand side remains bounded is therefore
$\frac1{1+3\gamma'}$.
Since $\gamma'>\gamma>3\frac{1+2\varepsilon}{1-2\varepsilon}$ are completely arbitrary,
this means that any value
$$\theta<\frac1{1+9\frac{1+2\varepsilon}{1-2\varepsilon}}=\frac{1-2\varepsilon}{10+16\varepsilon}$$
is achievable.
Letting $\varepsilon\downarrow 0$, any value
$\theta<\frac1{10}$ is achievable.

\nextpar{$H^{-1}$-estimate}
We have by (\ref{ineq:n-th-step-bound-on-H(-1)-norm})
\begin{eqnarray}
\sup_t\|v(\cdot,t)-v_0(\cdot,t)\|_{H^{-1}(\mathbb T^d)}
\leq r_0\sum_{n=0}^\infty\delta_n^\frac14
\leq\sigma\label{ineq:condition-on-r0-for-H(-1)-approximation}
\end{eqnarray}
by choosing $r_0$ sufficiently small.
\stopproof

\subsection{Proof of Corollary~\ref{corollary:strong-h-principle-for-2D-Euler}}
Corollary~\ref{corollary:strong-h-principle-for-2D-Euler} will follow from
\begin{proposition}\label{proposition:geometry-of-sym-pos-def-matrices}
With the notation of Section~\ref{section:notation}, 
if $d\geq 2$, then 
$$\mathcal M_d\subset \mathcal S^{d\times d}_{++}\qquad {\rm and}\qquad
{\rm Id}\in \mathcal M_d.$$
Furthermore, $R\in\mathcal S^{d\times d}_{++}$ is in $\mathcal M_d$ if and only if
$$\frac{{\rm tr}\,R}{d-1}\,{\rm Id}-R>0.$$
In particular,
$\mathcal M_2 =\mathcal S^{2\times 2}_{++}$ and 
$\mathcal M_d\subsetneq \mathcal S^{d\times d}_{++}$ for $d\geq 3$.
\end{proposition}
\startproof
It is obvious that $\mathcal M_d\subset \mathcal S^{d\times d}_{++}$
and one easily verifies that 
${\rm Id}=\frac1{d-1}\sum_{i=1}^d\left({\rm Id}-e_i\otimes e_i\right)$
where $\{e_1,\dots, e_d\}$ is the canonical basis for $\mathbb R^d$.
\bigskip

Suppose that $R\in \mathcal S^{d\times d}_{++}$ is of the form
$R=\sum_ia_i\left({\rm Id}-b_i\otimes b_i\right)$, where $a_i>0$ and $|b_i|=1$.
Then, $\sum_ia_i=\frac{{\rm tr}\,R}{d-1}$ and hence
$$0< \sum_ia_ib_i\otimes b_i=\frac{{\rm tr}\,R}{d-1}\,{\rm Id}-R.$$

Conversely, suppose $R\in\mathcal S^{d\times d}_{++}$ satisfies $\frac{{\rm tr}\,R}{d-1}\,{\rm Id}-R>0$.
$R$ is diagonalizable and 
all its eigenvalues satisfy $\lambda_i<\frac{{\rm tr}\,R}{d-1}$.
It is then easy to verify that, with $a_i:=\frac{{\rm tr}\,R}{d-1}-\lambda_i>0$, we have
after diagonalization
$$R=\sum_{i=1}^da_i\left({\rm Id}-e_i\otimes e_i\right)$$
where $\{e_1,\dots, e_d\}$ is the canonical basis of $\mathbb R^d$.
\bigskip

Finally, note that $\sum_{i=1}^d\lambda_i={\rm tr}\,R$ and $\lambda_i\geq 0$.
If $d=2$, then $\lambda_i\leq {\rm tr}\,R$ for $i=1,2$.
Otherwise, if $d\geq 3$, the condition $\frac{{\rm tr}\,R}{d-1}\,{\rm Id}-R>0$ can be violated.
\stopproof

\noindent{\it Proof of Corollary~\ref{corollary:strong-h-principle-for-2D-Euler}:}
It is easy to see that 
$R=\frac{e-\fint |v|^2}d\,{\rm Id}-\mathring R$ satisfies (\ref{ineq:strong-subsolution})
if and only if $R\in\mathcal M_d$.
Thus, Proposition~\ref{proposition:geometry-of-sym-pos-def-matrices} implies
Corollary~\ref{corollary:strong-h-principle-for-2D-Euler}.
\stopproof

\subsection{Proof of Corollary~\ref{corollary:3D-flows-are-not-2D}}
For sufficiently small $\sigma$ we have from the bound (\ref{ineq:int(vxv)})
$$|v_3|>
\frac12\fint_{\mathbb T^d}
\left(
|v_{0,33}(x,t)|^2+\frac{e(t)-\fint_{\mathbb T^d}|v_0(x',t)|^2}d-\mathring R_{0,33}(x,t)
\right)\,dx>\sigma
$$
for sufficiently small $\sigma$
whereas from (\ref{ineq:|v|H-1<epsilon}) we would have $|v_3|<\sigma$.
Thus, the solution cannot be of the form (\ref{eq:3D-flow-which-is-2D})
if $\sigma$ is chosen sufficiently small.

The analogous conclusion holds as well for the case $d=2$: the flows constructed
in Theorem~\ref{theorem:weak-h-principle} are genuinely two-dimensional,
that is, they are not parallel flows.
However, this conclusion can be arrived at by more elementary means.
Indeed, it is classical that such flows are necessarily stationary, 
and this is not possible if $e(t)$ is chosen non-constant.
\stopproof

\section{Construction of the iterates}
\subsection{Linear spaces of stationary flows}\label{section:linear-spaces-of-stationary-flows}
\setcounter{counterpar}{0}
An essential ingredient in the construction introduced in \cite{DS3}
is a \textit{linear} set of functions $(W,Q)$ (in the $\xi$-variable) 
satisfying the stationary Euler equations.
The existence of such spaces seems to hold for different reasons for $d=2$ and $d=3$,
and we consider these cases separately.

\nextpar{Dimension 3}
For $k\in\mathbb Z^3$, we let
\begin{equation}
b_k(\xi):=B_ke^{ik\cdot\xi},\qquad \psi_k(\xi):=D_ke^{ik\cdot\xi}
\label{def:bk(xi)-in-3D}
\end{equation}
where $B_k\in\mathbb C^3$ satisfies $|B_k|=\frac1{\sqrt 2}$, $k\cdot B_k=0$, and $\overline{B_k}=B_{-k}$,
and $D_k=i\frac{k\times B_k}{|k|^2}$ so that 
\begin{equation}
b_k={\rm curl}_\xi\,\psi_k,\qquad{\rm div}_\xi b_k=0,\qquad \overline{b_k}=b_{-k}.
\label{eq:properties-of-b(k)-in-3D}
\end{equation}
Here, the operator ${\rm curl}_\xi=\nabla_\xi\times \cdot$ is defined as usual.
The $\psi_k$ are \defn{vector potentials} for the vector fields $b_k$.
Concerning the analysis in this paper, they will play the same role as the stream functions in $d=2$ dimensions
introduced in (\ref{def:bk(xi)-in-2D}).

\nextpar{Dimension 2}
For $k\in \mathbb Z^2$, we let
\begin{equation}
b_k(\xi):=i\frac{k^\perp}{|k|}e^{ik\cdot\xi},
\qquad
\psi_k(\xi)=\frac{e^{ik\cdot\xi}}{|k|}
\label{def:bk(xi)-in-2D}
\end{equation}
so that
\begin{equation}
b_k(\xi)={\rm curl}_\xi\,\psi_k(\xi),\qquad {\rm div}_\xi\,b_k=0,\qquad \overline{b_k}=b_{-k}
\label{eq:properties-of-b(k)-in-2D}
\end{equation}
where this time
${\rm curl}_\xi=\nabla_\xi^\perp=(-\partial_{\xi^2},\partial_{\xi^1})$ denotes the rotated gradient.
From the analytic point of view, 
the \defn{stream function} $\psi_k$ is the analogue of the vector potential $D_ke^{ik\cdot\xi}$ defined 
in (\ref{def:bk(xi)-in-3D})
for the case $d=3$.
\bigskip

\begin{lemma}\label{lemma:linear-spaces-of-steady-flows}
Let $\nu\geq 1$ and $d=2$ or $3$.
For $k\in\mathbb Z^d$ such that $|k|^2=\nu$, let $a_k\in\mathbb C$ such that $\overline{a_k}=a_{-k}$.
Then
\begin{equation}
W(\xi)=\sum_{|k|^2=\nu}a_kb_k(\xi),
\qquad Q
:=\left\{\begin{array}{ll}-\frac{|W|^2}2+\fint \frac{|W|^2}2\,d\xi&(d=3)\\
-\frac{|W|^2}2+\nu\frac{\Psi^2}2&(d=2)
\end{array}
\right.,
\label{def:W-bk}
\end{equation}
where
$\Psi(\xi)=\sum_{|k|^2=\nu}a_k\psi_k(\xi)$,
are $\mathbb R$-valued and 
satisfy
\begin{equation}
{\rm div}_\xi\,\left(W\otimes W\right)+\nabla_\xi\,Q=0,\quad {\rm div}_\xi\,W=0.
\label{eq:div(WxW)+nabla(Q)=0}
\end{equation}
Furthermore,
\begin{equation}
\fint_{\mathbb T^d}W\otimes W\,d\xi=\sum_{|k|^2=\nu}|a_k|^2\left({\rm Id}-\frac k{|k|}\otimes \frac k{|k|}\right).
\label{eq:ave(WxW)=sum|ak|2Mk}
\end{equation}
\end{lemma}
\startproof
If $d=3$, this is Lemma~3.1 of \cite{DS3}.
(A constant is added in our definition so that $\fint_{\mathbb T^3} Q\,d\xi=0$.)
Suppose $d=2$.
By direct computation one finds $\Delta_\xi\psi_k=-|k|^2\psi_k$, 
and hence that $\Delta_\xi \Psi=-\nu\Psi$.
Recall the identities
$${\rm div}\,_\xi(W\otimes W)=\frac12\nabla_\xi|W|^2+({\rm curl}_\xi W)W^\perp$$
where ${\rm curl}_\xi W=\partial_{\xi^1}W^2-\partial_{\xi^2}W^1=\Delta_\xi\Psi$ and $W^\perp=(-W^2,W^1)$.
Then,
$${\rm div}\,_\xi(W\otimes W)=\nabla_\xi\frac{|W|^2}2-\nu\Psi\nabla_\xi\Psi$$
as desired.

As for the average, write
\begin{eqnarray}
\fint_{\mathbb T^2}W\otimes  W(\xi)\,d\xi
=\sum_{j,k}a_k\overline{a_j}\fint_{\mathbb T^d}e^ {i(k-j)\cdot \xi}\,d\xi\frac{j^\perp}{|j|}\otimes \frac{k^\perp}{|k|}
=\sum_{|k|^2=\nu}|a_k|^2\left({\rm Id}-\frac k{|k|}\otimes \frac k{|k|}\right)
\notag
\end{eqnarray}
where the last identity follows from $\frac{k^\perp}{|k|}\otimes \frac{k^\perp}{|k|}=\left({\rm Id}-\frac k{|k|}\otimes \frac k{|k|}\right)$
by direct calculation.
\stopproof

\subsection{The Geometric Lemma}\label{section:geometric-lemma}
The next Lemma is a quantified examination of the range of positive definite matrices
that the flows from Lemma~\ref{lemma:linear-spaces-of-steady-flows} are able to generate.

\begin{lemma}[Geometric Lemma]\label{lemma:geometric-lemma}
Suppose $d\geq 2$ and $N\geq 1$.
Let $K\subset \mathcal M_d$ be compact, and $\mathcal N$ an open neighborhood of $K$
such that $\overline{\mathcal N}\subset \mathcal M_d$.
Then, there exist $\nu\geq 1$, 
pairwise disjoint subsets
$$
\Lambda_j\subset 
\{k\in \mathbb Z^d~:~|k|^2=\nu\}
\qquad j\in\{1,\dots, N\}
$$
and smooth positive functions
$$
\gamma_k^{(j)}\in C^\infty(\mathcal N),\qquad j\in \{1, \dots, N\}, \quad k\in \Lambda_j
$$
such that
\begin{enumerate}
\item $k\in\Lambda_j$ implies $-k\in\Lambda_j$ and $\gamma_k^{(j)}=\gamma_{-k}^{(j)}$;
\item for each $R\in \mathcal N$ and $j=1,2,\dots, N$ we have
\begin{equation}
R=\sum_{k\in\Lambda_j}\left(\gamma_k^{(j)}(R)\right)^2\left({\rm Id}-\frac k{|k|}\otimes \frac k{|k|}\right).
\label{eq:R=sum-of-Mk}
\end{equation}
\end{enumerate}
\end{lemma}
\startproof
Each $R\in\overline{\mathcal N}\subset \mathcal M_d$ is in the interior of a simplex $\Sigma(R)$ 
with vertices of the form
$$A_i(R)=\kappa_i(R)\left({\rm Id}-b_i(R)\otimes b_i(R)\right), \qquad i=1,\dots, d+1$$
where $|b_i(R)|=1$ and $\kappa_i(R)>0$.
$\Sigma(R)$ can be chosen with pairwise distinct vertices and hence $R$ is of the form
$$R=\sum_{i=1}^{d+1}c_i(R)\left({\rm Id}-b_i(R)\otimes b_i(R)\right)$$
where $c_i(R)>0$.

Since $\overline{\mathcal N}$ is compact,
we may extract a finite subcover $\{\Sigma_l\}_{l=1}^L$ where $\Sigma_l:=\Sigma(R_l)$.
Observe now that since $R_l$ is in the interior of $\Sigma_l$, it is also 
in the interior of any simplex with vertices slightly perturbed.
Recall now that $\mathbb Q^d\cap \mathbb S^{d-1}$ is dense in $\mathbb S^{d-1}$
(the proof in \cite{DS3} using stereographic projection holds in any dimension).
Then, by taking $\nu\in\mathbb N$ sufficiently large,
there exist $k_{i,l}^{(j)}\in\mathbb Z^d$, $i=1,\dots, d+1$, $l=1,\dots,L$, $j=1,\dots, N$, all distinct, satisfying $|k_{i,l}^ {(j)}|^2=\nu$,
and such that, for each $l=1,\dots, L$, $R_l$ is in the interior of $\Sigma_l^{(j)}$
for $j=1,\dots, N$, where $\Sigma_l^{(j)}$ is a simplex 
whose vertices are multiples of
$\left({\rm Id}-\frac{k_{i,l}^{(j)}}{|k_{i,l}^{(j)}|}\otimes \frac{k_{i,l}^{(j)}}{|k_{i,l}^{(j)}|}\right)$.
We then write
$$R_l=\sum_{i=1}^{d+1}c_{i,l}^{(j)} \left({\rm Id}-\frac{k_{i,l}^{(j)}}{|k_{i,l}^{(j)}|}\otimes \frac{k_{i,l}^{(j)}}{|k_{i,l}^{(j)}|}\right),
\qquad l=1,\dots, L,\quad j=1,\dots, N
$$
where $c_{i,l}^{(j)}>0$, $l=1,\dots,L$, $i=1,\dots,d+1$.

For each $l=1, \dots, L$, $j=1,\dots, N$
there exist positive functions $\alpha_{i,l}^{(j)}\in C^\infty(\Sigma_l)$, $i=1,\dots, d+1$, such that
$$R=
\sum_{i=1}^{d+1}\alpha_{i,l}^{(j)}(R)\left({\rm Id}-\frac{k_{i,l}^{(j)}}{|k_{i,l}^{(j)}|}\otimes \frac{k_{i,l}^{(j)}}{|k_{i,l}^{(j)}|}\right)
\qquad {\rm for}\quad R\in\Sigma_l^{(j)}.$$
(Indeed, $R\in\Sigma_l^{(j)}$ is the unique convex combination of the vertices of $\Sigma_l^{(j)}$,
and the coefficients are algebraic expressions of $R$.)
For each $j=1,\dots, N$, let now $\{\eta_l^{(j)}\}_{l=1}^L$ be a $C^\infty$ partition of unity subordinate to the cover $\{\Sigma_l^{(j)}\}_{l=1}^L$.
Then for any $R\in\overline{\mathcal N}$ we have
$$R=\sum_{l=1}^L\sum_{i=1}^{d+1}\eta_l^{(j)}(R)\alpha_{i,l}^{(j)}(R)\left({\rm Id}-\frac{k_{i,l}^{(j)}}{|k_{i,l}^{(j)}|}\otimes \frac{k_{i,l}^{(j)}}{|k_{i,l}^{(j)}|}\right),\qquad j=1,\dots,N.$$

Set now
$$\Gamma_j:=\left\{k_{i,l}^{(j)}\in\mathbb Z^d~\colon~ i=1,\dots,d+1, l=1,\dots,N\right\},\qquad j=1,\dots, N$$
and
$a_k:=\sqrt{\eta_l^{(j)}\alpha_{i,l}^{(j)}}$ for $k=k_{i,l}^{(j)}$.
Next let
$$\Lambda_j:=\Gamma_j\cup (-\Gamma_j),\qquad j=1,\dots, N$$
(once again by density, we can arrange for the $\Lambda_j$'s to be pairwise disjoint)
and $a_k=0$ if $k\not\in\Gamma_j$.
Finally, taking
$$\gamma_k(R):=\frac{\sqrt 2}2\sqrt{a_k(R)+a_{-k}(R)},\qquad k\in\Lambda_j$$
finishes the proof.
\stopproof

\subsection{The operator $\mathcal R={\rm div}^{-1}$}
\begin{definition}[The Leray projector]\label{definition:Leray-projectors}
Let $d\geq 2$.
For a vector field $v\in C^\infty(\mathbb T^d,\mathbb R^d)$, 
set
$$
\mathcal Qv:=\nabla \phi+\fint_{\mathbb T^d}v
$$
where $\phi\in C^\infty(\mathbb T^d)$ is the solution to
$\Delta\phi={\rm div}\,v$ in $\mathbb T^d$ subject to $\fint_{\mathbb T^d}\phi=0$.
We denote by
$\mathcal P:=I-\mathcal Q$
the Leray projector onto divergence-free vector fields with zero average.
\end{definition}
The operator $\mathcal R$ was introduced in \cite{DS3} for $d=3$.
Its generalization for any $d\geq 2$ is given by the following
\begin{definition}[The operator $\mathcal R$]\label{definition:the-operator-R}
Let $d\geq 2$.
For any smooth vector field $v\in C^\infty(\mathbb T^d,\mathbb R^d)$,
we define $\mathcal Rv$ to be the matrix-valued periodic function
\begin{equation}
\mathcal Rv
=\frac{d-2}{2(d-1)}\left(\nabla\mathcal Pu +\left(\nabla\mathcal P u\right)^\top\right)
+\frac d{2(d-1)}\left(\nabla u+(\nabla u)^\top\right)+\frac1{1-d}({\rm div}u){\rm Id}
\label{eq:R=divInverse}
\end{equation}
where $u\in C^\infty(\mathbb T^d,\mathbb R^d)$ is the solution to
$$
\Delta u=v-\fint_{\mathbb T^d}v\quad {\rm in}~\mathbb T^d,\qquad {\rm subject~to}\quad \fint_{\mathbb T^d}u=0.
$$
\end{definition}

By direct verification one obtains
\begin{lemma}[$\mathcal R={\rm div}\,^{-1}$]\label{lemma:R=divInverse}
Let $d\geq 2$.
For any $v\in C^\infty(\mathbb T^d,\mathbb R^d)$ we have
\begin{enumerate}
\item $\mathcal Rv(x)$ is a symmetric trace-free matrix for each $x\in\mathbb T^d$;
\item ${\rm div}\,\mathcal Rv=v-\fint_{\mathbb T^d}v$.
\end{enumerate}
\end{lemma}

\subsection{Further technical preliminary}
The following is proved in \cite{DS3} for $d=3$ and the proof is valid as it is for any number of space dimensions.
Denote $\mathcal C_1, \dots, \mathcal C_{2^d}$ the equivalence classes 
of $\mathbb Z^d/\sim$ where $k\sim l$ if $k-l\in (2\mathbb Z)^d$.
\begin{proposition}[Partition of the space of velocities]
Let $d\geq 2$ and $\mu\in\mathbb N$.
There exists a partition of the space of velocities, 
namely $\mathbb R$-valued functions $\alpha_l(v)$ for $l\in\mathbb Z^d$ satisfying
\begin{equation}
\sum_{l\in\mathbb Z^d}(\alpha_l(v))^2\equiv 1
\label{eq:sum-of-ak-squared=1}
\end{equation}
such that, setting
$\phi_k^{(j)}(v,\tau)=\sum_{l\in\mathcal C_j}\alpha_l(\mu v)e^{-i(k\cdot\frac l\mu)\tau}$,
for 
$j=1,\dots, 2^d$, and $k\in\mathbb Z^d$,
then we have $\overline{\phi_k^{(j)}}=\phi_{-k}^{(j)}$ and
\begin{equation}|\phi_k^{(j)}(v,\tau)|^2=\sum_{l\in\mathcal C_j}\alpha_l^2(v).
\label{eq:|phi(j)k|2=sum-of-alpha(k)2}
\end{equation}
\end{proposition}
\bigskip

\subsection{The maps $w_o$, $v_1$, $p_1$, and $\mathring R_1$}
\label{section:definitions-of-the-maps}
\setcounter{counterpar}{0}
Let $e_0,\Delta_0,\zeta, (v,\mathring R,p), r_0, K, \mathcal N, \delta$
be as in Proposition~\ref{proposition:iteration}.

\nextpar{Mollifications}
Let $\chi\in C^\infty_c(\mathbb R^d\times \mathbb R)$ be a smooth standard nonnegative radial kernel supported in 
$\left[-\pi,\pi\right]^{d+1}$ and denote by
$$\chi_\ell(x,t):=\frac1{\ell^{d+1}}\chi(\frac x\ell,\frac t\ell)$$
the corresponding family of mollifiers ($0<\ell<1$).
We define
\begin{equation}
v_\ell(x,t):=\int_{\mathbb T^d\times \mathbb S^1}v(x-y,t-s)\chi_\ell(y,s)\,dy\,ds,
\label{eq:v(l)}
\end{equation}
$$\mathring R_\ell(x,t):=\int_{\mathbb T^d\times \mathbb S^1}\mathring R(x-y,t-s)\chi_\ell(y,s)\,dy\,ds$$
and set
\begin{equation}
\rho_\ell(t) := \frac1{d(2\pi)^d}\left(e_0(t)+\Delta_0(t)(1-\overline\delta)-\fint_{\mathbb T^d}|v_\ell(x,t)|^2\,dx
\right)
\label{eq:rho(l)}
\end{equation}
and
\begin{equation}R_\ell(x,t):=\rho_\ell(t){\rm Id}-\mathring R_\ell(x,t),\qquad x\in\mathbb T^d,\quad t\in\mathbb S^1.
\label{eq:R(l)}
\end{equation}
\bigskip

\nextpar{The oscillation term $w_o$}
Provided $\frac{R_\ell}{\rho_\ell}\in\mathcal N$, see Lemma~\ref{lemma:wo-is-well-defined},
we may define
\begin{equation}
w_o(x,t):=W(x,t;\lambda t, \lambda x)
\label{eq:wo(x,t)}
\end{equation}
where (the $b_k$'s are defined in (\ref{def:bk(xi)-in-2D}) and (\ref{def:bk(xi)-in-3D}))
\begin{eqnarray}
W(y,s;\tau,\xi)
&:=&\sum_{|k|^2=\nu}a_k(y,s;\tau)b_k(\xi)\notag\\
&:=&\sqrt{\rho_\ell(s)}\sum_{j=1}^{2^d}\sum_{k\in\Lambda_j}\gamma_k\left(\frac{R_\ell(y,s)}{\rho_\ell(s)}\right)\phi_{k,\mu}^{(j)}(v_\ell(y,s),\tau)b_k(\xi).
\label{eq:W(y,s,tau,xi)}
\end{eqnarray}
and where for $k\in\Lambda_j$,
\begin{equation}
a_k(y,s;\tau)=\sqrt{\rho_\ell(s)}
\gamma_k\left(\frac{R_\ell(y,s)}{\rho_\ell(y,s)}\right) \phi_{k,\mu}^{(j)}(v_\ell(y,s);\tau)
\label{eq:ak}
\end{equation}

The corresponding stream function ($d=2$) and vector potential ($d=3$) are both formally defined by
\begin{equation}
\psi_o(x,t):=\Psi(x,t;\lambda t, \lambda x)
\label{eq:psio(x,t)}
\end{equation}
where (see again (\ref{def:bk(xi)-in-2D}) and (\ref{def:bk(xi)-in-3D}))
\begin{eqnarray}
\Psi(y,s;\tau,\xi)
&:=&\sum_{|k|^2=\nu}a_k(y,s;\tau)\psi_k(\xi)
\notag.
\end{eqnarray}

\nextpar{The velocity field $v_1$}
It is defined by
$$v_1:=v+w:=v+w_o+w_c,\qquad w_c:=-\mathcal Q w_o$$
where $w_o$ is given in (\ref{eq:wo(x,t)}) and $\mathcal Q$ is the Leray projector of
Definition~\ref{definition:Leray-projectors}.
Note that ${\rm div}\,v_1=0$.

\nextpar{The pressure $p_1$}
It is defined by
\begin{equation}
p_1:=p+q=p+\tilde p-2\frac{(v-v_\ell)\cdot w}d
\label{eq:p1}
\end{equation}
where
\begin{equation}
\tilde q(x,t):=Q(x,t;\lambda t,\lambda x)
\label{eq:q(tilde)}
\end{equation}
and
\begin{equation}
Q(y,s;\tau,\xi)
:=\sum_{1\leq |k|\leq 2\nu}\tilde a_k(y,s;\tau)e^{ik\cdot\xi}
:=\left\{\begin{array}{ll}-\frac{|W|^2}2+\fint \frac{|W|^2}2\,d\xi&(d=3)\\
-\frac{|W|^2}2+\nu\frac{\psi_o^2}2&(d=2)
\end{array}
\right.
\label{eq:Q}
\end{equation}
and where $v_\ell$ is given in (\ref{eq:v(l)}), $w_o$ in (\ref{eq:wo(x,t)}), 
$W$ in (\ref{eq:W(y,s,tau,xi)}),
and in case $d=2$, $\nu$ is given by Geometric Lemma~\ref{lemma:geometric-lemma}.

\nextpar{The tensor $\mathring R_1$}
We define
\begin{eqnarray}\mathring R_1
&:=& \mathring R-\mathring R_\ell\notag\\
&&+w\otimes (v-v_\ell)+(v-v_\ell)\otimes w-2\frac{(v-v_\ell)\cdot w}d{\rm Id}\notag\\
&&+\mathcal R\,[{\rm div}\,(w_o\otimes w_o+\mathring R_\ell+\tilde q{\rm Id})\notag\\
&&+\mathcal R\,\partial_tw_c\notag\\
&&+\mathcal R\,{\rm div}\,\left((v_\ell+w)\otimes w_c+w_c\otimes (v_\ell+w)-w_c\otimes w_c\right)\notag\\
&&+\mathcal R\,{\rm div}\,(w_o\otimes v_\ell)\notag\\
&&+\mathcal R\,\left[\partial_tw_o+{\rm div}\,(v_\ell\otimes w_o)\right]
\quad \left(=\mathcal R\left[\partial_tw_o+v_\ell\cdot \nabla w_o\right]\right).\notag\\
&&\label{eq:Ro1}
\end{eqnarray}
One easily verifies (see \S~3.5 in \cite{DS-Hoelder} for details)
that $\mathring R_1\in\mathcal S^{d\times d}_0$ and that 
$${\rm div}\,\mathring R_1=\partial_tv_1+{\rm div}\,(v_1\otimes v_1)+\nabla p_1.$$

\section{Proof of Proposition~\ref{proposition:iteration}}\label{section:proof-of-Proposition(iteration)}
\subsection{Conditions on the parameters}
Let $e_0,\Delta_0, \overline r_0, \varepsilon, \zeta, K,\mathcal N$ 
be as in Proposition~\ref{proposition:iteration}.
Set $e:=e_0+\Delta_0$ and $\omega:=\frac \varepsilon{2+\varepsilon}$ so that
$$1+\varepsilon=\frac{1+\omega}{1-\omega}.$$
We assume $D\geq 1$ and $\delta\leq 1$ are given.

The estimates in the following Section as well as those established in \cite{DS-Hoelder},
see Propostions~\ref{proposition:estimates-on-wo-wc-w}, \ref{proposition:estimates-on-the-energy},
and \ref{proposition:estimates-on-Ro1} in the Appendix,
are derived under the assumptions on the parameters $\lambda, \mu$ and $\ell$
that they satisfy
\begin{equation}\lambda,\mu,\frac\lambda\mu\in\mathbb N
\notag
\end{equation}
and
\begin{equation}
\mu\geq \delta^{-1}\geq 1,\qquad \ell^{-1}\geq \frac D{\eta\delta}\geq 1,
\qquad \lambda \geq \max\left\{(\mu D)^{1+\omega},\ell^{-(1+\omega)}\right\}.
\label{ineq:conditions-on-the-parameters}
\end{equation}

\subsection{$w_o$ is well defined}
The following estimates are standard:
\begin{eqnarray}
\|v_\ell\|_r&\leq& C(r)D\ell^{-r}\qquad (r\geq 1),
\label{ineq:estimates-on-v(l)}\\
\|v-v_\ell\|_0+\|\mathring R_\ell-\mathring R\|_0&\leq&  CD\ell,
\label{ineq:estimates-on-|v-v(l)|-and-|Ro-Ro(l)|}\\
\|\mathring R_\ell\|_0&\leq&\|\mathring R\|_0.
\label{ineq:|Ro(l)|<|Ro|}
\end{eqnarray}
As a consequence, writing $\left||v_\ell|^2-|v|^2\right|\leq |v-v_\ell|^2+2|v||v-v_\ell|$,
and using $D\ell\leq \eta\delta\leq \frac{\Delta_0}{4d}r_0$ from (\ref{ineq:conditions-on-the-parameters})
and (\ref{eq:eta}) 
we obtain
\begin{equation}
\fint_{\mathbb T^d}\left||v_\ell(x,t)|^2-|v(x,t)|^2\right|\,dx
\leq CD\ell \left(\max e^\frac12+1\right)
\leq C\eta\delta \left(\max e^\frac12+1\right).
\label{ineq:int||v(l)|2-|v|2|}
\end{equation}

\begin{lemma}[$w_o$ is well defined]\label{lemma:wo-is-well-defined}
Let $d, e_0, \Delta_0, K, \mathcal N, \overline r_0$ as in Proposition~\ref{proposition:iteration}.
Let $0<\delta\leq 1$, $\overline \delta\leq \zeta\delta$, $0<\zeta\leq \frac12$,
and $D\ell\leq\delta\frac{\min\Delta_0}{4d}\overline r_0$.
Then there exists $\tilde r_0$ depending on $d, K, e_0$, 
and $\Delta_0$ from Proposition~\ref{proposition:iteration},
such that the following holds.
If $r_0\leq \tilde r_0$,
$(v, \mathring R,p)$ satisfies (\ref{ineq:assumed-energy-gap})
and (\ref{eq:rescaled-R-is-in-K}),
and if $\rho_\ell$ and $R_\ell$ are defined as in (\ref{eq:rho(l)}) and (\ref{eq:R(l)}) respectively,
then
$\frac{R_\ell}{\rho_\ell}\in \overline{\mathcal N}$.
\end{lemma}
\startproof
By assumption
${\rm Id}-\frac d{e_0(t)+\Delta_0(t)(1-\overline\delta)-\fint_{\mathbb T^d}|v(x,t)|^2\,dx}\mathring R\in K$.
In order to prove that
${\rm Id}-\frac d{e_0(t)+\Delta_0(t)(1-\overline\delta)-\fint_{\mathbb T^d}|v_\ell(x,t)|^2\,dx}\mathring R_\ell\in \mathcal N$,
we shall prove that
\begin{eqnarray}
\left\|\frac{d\mathring R}{e_0(t)+\Delta_0(t)(1-\overline\delta)-\fint_{\mathbb T^d}|v|^2}
-\frac{d\mathring R_\ell}{e_0(t)+\Delta_0(t)(1-\overline\delta)-\fint_{\mathbb T^d}|v_\ell|^2}\right\|
\notag
\end{eqnarray}
is less than
${\rm dist}\,(K,\partial\mathcal N):=\inf\left\{|A-B|~\colon~A\in K, B\in\partial\mathcal N\right\}$.
By assumptions on $\overline\delta$ and $\zeta$,
$$e_0+\Delta_0(1-\overline\delta)-\fint_{\mathbb T^d}|v|^2\geq \Delta_0\delta/4$$
and thus
\begin{eqnarray}
&&d(2\pi)^d\rho_\ell(t)\notag\\
&=& e_0(t)+\Delta_0(t)(1-\overline\delta)-\fint_{\mathbb T^d}|v_\ell(x,t)|^2\,dx\notag\\
&=& e_0(t)+\Delta_0(t)(1-\overline\delta)-\fint_{\mathbb T^d}|v(x,t)|^2\,dx
-\fint_{\mathbb T^d}\left(|v_\ell(x,t)|^2-|v(x,t)|^2\right)\,dx\notag\\
&\geq&\Delta_0(t)\delta/4-\fint_{\mathbb T^d}\left||v_\ell(x,t)|^2-|v(x,t)|^2\right|\,dx.
\label{ineq:condition-on-eta-for-lower-bound-on-rho(l)}
\end{eqnarray}
By (\ref{ineq:int||v(l)|2-|v|2|}),
making $\overline r_0$ smaller if necessary depending on $d$, $\Delta_0$ and $e=e_0+\Delta_0$,  we have
\begin{eqnarray}
d(2\pi)^d\rho_\ell(t)\geq \Delta_0(t)\delta/8.
\notag
\end{eqnarray}

Since $(v,\mathring R,p)$ satisfies (\ref{eq:rescaled-R-is-in-K}),
there exists a constant $C=C(K)$ such that
$$\left\|\frac{\mathring R}{e_0+\Delta(1-\overline\delta)-\fint_{\mathbb T^d}|v|^2}\right\|\leq C(K).$$

Using the above, (\ref{ineq:estimates-on-|v-v(l)|-and-|Ro-Ro(l)|}), (\ref{ineq:int||v(l)|2-|v|2|}), 
and again $D\ell\leq \frac{\min\Delta_0}d\overline r_0$, we obtain
\begin{eqnarray}
&&\left|\frac{\mathring R}{e_0+\Delta_0(1-\overline\delta)-\fint_{\mathbb T^d}|v|^2}
-\frac{\mathring R_\ell}{e_0+\Delta_0(1-\overline\delta)-\fint_{\mathbb T^d}|v_\ell|^2}
\right|\notag\\
&\leq&\left|\frac{\mathring R-\mathring R_\ell}{e_0+\Delta_0(1-\overline\delta)-\fint_{\mathbb T^d}|v|^2}\right|\notag\\
&&+|\mathring R_\ell|\left|\frac1{e_0+\Delta(1-\overline\delta)-\fint_{\mathbb T^d}|v|^2}-\frac1{e_0+\Delta(1-\overline\delta)-\fint_{\mathbb T^d}|v_\ell|^2}\right|
\notag\\
&\leq&\left|\frac{\mathring R-\mathring R_\ell}{e_0+\Delta_0(1-\overline\delta)-\fint_{\mathbb T^d}|v|^2}\right|\notag\\
&&+\frac{|\mathring R|}{e_0+\Delta_0(1-\overline\delta)-\fint_{\mathbb T^d}|v|^2}\;
\frac{\fint_{\mathbb T^d}\left||v_\ell|^2-|v|^2\right|}{e_0+\Delta_0(1-\overline\delta)-\fint_{\mathbb T^d}|v_\ell|^2}\notag\\
&\leq&\frac{4CD\ell}{\Delta_0\delta}
+C(K)\frac{4C\eta\delta (\max e^\frac12+1)}{\Delta_0\delta}\notag\\
&\leq&\frac{\overline r_0}d\left(C+C(K)(\max e^\frac12+1)\right)
\label{ineq:condition-on-eta-so-that-R(l)/rho(l)-is-in-N}
\end{eqnarray}
Therefore, the right-hand side is sufficiently small so that
$\frac{R_\ell}{\rho_\ell}\in\mathcal N$,
provided $\overline r_0\leq \tilde r_0$ where $\tilde r_0$ is chosen sufficiently small depending on $d, K, e, \Delta_0$.
\stopproof

\subsection{Proof of Proposition~\ref{proposition:iteration}}
\setcounter{counterpar}{0}
\smallskip

\nextpar{Setting some parameters}
In the next paragraphs, we will use estimates from \cite{DS-Hoelder}, 
see Propositions~\ref{proposition:estimates-on-wo-wc-w}, \ref{proposition:estimates-on-the-energy},
and \ref{proposition:estimates-on-Ro1} in the Appendix.
These estimates are derived under the conditions listed in 
(\ref{ineq:conditions-on-the-parameters}) on the parameters $\ell,\lambda,\mu,D$, and $\varepsilon$ ({\it via} $\omega$).
We shall now {\it set} the parameters $\ell, \mu, \lambda$ in terms of $D, \delta, \varepsilon$
so that these conditions are satisfied.
Set
$$\alpha=\frac\omega{2(1+\omega)}.$$
In particular, both $\omega$ and $\alpha$ depend only on $\varepsilon$
and $\alpha\in(0,\frac\omega{1+\omega})$ so that Propositions~\ref{proposition:estimates-on-the-energy}
and \ref{proposition:estimates-on-Ro1} are applicable.
Note then that the constants $C_{v,s}$ become constants $C_v$.
Also,
$$\alpha-\frac12=-\frac1{2(1+\omega)}\quad<\quad 2\alpha-\frac12=-\frac{1-\omega}{2(1+\omega)}\quad<\quad 0.$$

Recall that $\overline\delta=\zeta\delta^\frac32$ and choose
\begin{equation}\ell=\frac1{L_v}\frac{\overline\delta}D
\label{eq:ell}
\end{equation}
where $L_v\geq 1$ will be chosen sufficiently large, see Section~\ref{section:fixing-the-parameters}.
We shall impose
\begin{equation}\mu^2D=\lambda
=\Lambda_v\left(\frac{D\delta}{\overline\delta^2}\right)^\frac1{1-4\alpha}
=\Lambda_v\left(\frac{D\delta}{\overline\delta^2}\right)^\frac{1+\omega}{1-\omega}
=\Lambda_v\left(\frac{D\delta}{\overline\delta^2}\right)^{1+\varepsilon}
\label{eq:simplifying-relation-between-parameters}
\end{equation}
where $\Lambda_v\geq 1$ will be chosen sufficiently large, see Section~\ref{section:fixing-the-parameters}.
(We note that in principle we should require that $\lambda, \mu,\lambda/\mu\in\mathbb N$,
but this can be arranged easily, up to universal constants.)


Now we verify that the conditions (\ref{ineq:conditions-on-the-parameters}) on the parameters
are satisfied with the above choices (\ref{eq:ell}) and (\ref{eq:simplifying-relation-between-parameters}).
Noting that $\overline\delta\leq \delta$, then
$\ell^{-1}\geq \frac D{\eta\delta}$ is satisfied with
\begin{equation}L_v\geq \eta^{-1}.
\label{ineq:condition-on-Lv-to-achieve-conditions-on-parameters}
\end{equation}

Next, (\ref{eq:simplifying-relation-between-parameters}) and $\Lambda_v\geq 1$ imply
$$\mu=\sqrt{\frac\mu{D}}=\Lambda_v^\frac12\left(\frac D{\zeta^2\delta^2}\right)^\frac{1+\varepsilon}2D^{-\frac12}
=\frac{\Lambda_v^\frac12D^\frac\varepsilon2}{\zeta^{1+\varepsilon}}\delta^{-(1+\varepsilon)}\geq \delta^{-1}
$$
since $\zeta\leq \frac12$, and $D\geq 1$.
Also,
$$\frac\lambda{(\mu D)^{1+\omega}}=\lambda^\frac{1-\omega}2D^\frac{1+\omega}2
=\Lambda_v^\frac{1-\omega}2\left(\frac D{\zeta^2\delta^2}\right)^\frac{1+\omega}2D^\frac{1+\omega}2
=\frac{\Lambda_v^\frac{1-\omega}2}{\zeta^{1+\omega}}D^{1+\omega}\delta^{-(1+\omega)}\geq 1
$$
since $0<\omega<1$, $\delta\leq 1$, $\Lambda_v\geq 1$, and $D\geq 1$.
Also, 
\begin{eqnarray}
=\Lambda_v
\left(\frac{D}{\zeta^2\delta^2}\right)^{\frac{1+\omega}{1-\omega}}
\left(\frac{\zeta\delta^\frac32}{L_vD}\right)^{\frac{1+\omega}{1-\omega}}
=\frac{\Lambda_v}{\zeta^{\frac{(1+\omega)^2}{1-\omega}}L_v^{1+\omega}} 
\geq
\frac{\Lambda_v}{\zeta^{\frac{(1+\omega)^2}{1-\omega}}L_v^{1+\omega}}
\notag
\end{eqnarray}
so that we require
\begin{equation}
\frac{\Lambda_v}{\zeta^{\frac{(1+\omega)^2}{1-\omega}}L_v^{1+\omega}}\geq 1
\label{ineq:condition-on-Lv-Lambdav-to-achieve-conditions-on-parameters}
\end{equation}

In conclusion, the requirements (\ref{ineq:conditions-on-the-parameters}) are satisfied
provided $L_v\geq 1$ and $\Lambda_v\geq 1$ satisfy (\ref{ineq:condition-on-Lv-to-achieve-conditions-on-parameters}),
and (\ref{ineq:condition-on-Lv-Lambdav-to-achieve-conditions-on-parameters}).
Note that $\zeta$ shall be chosen first, then $L_v$, and finally $\Lambda_v$.
(Further requirements will be imposed on $\zeta$, $L_v$ and $\Lambda_v$,
see Section~\ref{section:fixing-the-parameters}.)

\nextpar{Estimates on the energy}
From Proposition~\ref{proposition:estimates-on-the-energy}, and with
$\alpha=\frac\omega{1+\omega}$,
\begin{eqnarray}
&&\left|e_0(t)+\Delta_0(1-\overline \delta)-\fint_{\mathbb T^d}|v_1(x,t)|^2\,dx\right|\notag\\
&\leq&C_eD\ell+C_v\sqrt\delta D^\frac12\lambda^{\alpha-\frac12}\notag\\
&\leq&\frac{C_e}{L_v}\overline\delta
+\frac{C_v}{\Lambda_v^\frac1{2(1+\omega)}}D^{\frac12-\frac1{2(1-\omega)}}\delta^{\frac12-\frac1{2(1-\omega)}}\overline\delta^\frac1{1-\omega}\notag\\
&\leq&\frac{C_e}{L_v}\overline\delta
+\frac{C_v}{\Lambda_v^\frac1{2(1+\omega)}}\overline\delta
\notag
\end{eqnarray}
where simplifications follow since $D\geq 1$, $\frac12-\frac1{2(1-\omega)}=-\frac\omega{2(1-\omega)}<0$ for $0<\omega<1$,
and
$$\delta^{\frac12-\frac1{2(1-\omega)}}\overline\delta^\frac1{1-\omega}=\delta^{-\frac\omega{2(1-\omega)}}\overline\delta^\frac\omega{1-\omega}\overline\delta
=\left(\frac{\overline\delta}{\delta^\frac12}\right)^\frac\omega{1-\omega}\overline\delta\leq \overline\delta$$
since $\overline\delta\leq \sqrt\delta$.
We can achieve (\ref{ineq:target-energy-gap}) provided
\begin{equation}
\frac{C_e}{L_v}+\frac{C_v}{\Lambda_v^\frac1{2(1+\omega)}}
\leq \zeta/2\min_t\Delta_0(t).
\label{ineq:condition-on-Lv-Lambdav-to-achieve-target-energy-gap}
\end{equation}

Using $\frac{\alpha-\frac12}{1-4\alpha}=-\frac1{2(1-\omega)}$
and $\delta^{\frac12-\frac1{2(1-\omega)}}\overline\delta^\frac1{1-\omega}\leq \overline\delta$, established above,
the second estimate in Proposition~\ref{proposition:estimates-on-the-energy} becomes
\begin{eqnarray}
&&\left|\int_{\mathbb T^d}\left(v_1\otimes v_1-v\otimes v-R_\ell\right)\,dx\right|\notag\\
&\leq&
C_{v,e}\sqrt\delta D^\frac12\lambda^{\alpha-\frac12} +C_e\delta D^\frac12\lambda^{-\frac12}\notag\\
&\leq&
C_{v,e}\Lambda_v^{-\frac1{2(1+\omega)}} D^{\frac12{-\frac1{2(1-\omega)}}}\delta^{\frac12-\frac1{2(1-\omega)}}\overline\delta^\frac1{1-\omega}
+C_e\Lambda_v^{-\frac12}\delta^\frac12\overline\delta\notag\\
&\leq&\left(\frac{C_{v,e}}{\Lambda_v^{\frac1{2(1+\omega)}}}+\frac{C_e}{\Lambda_v^\frac12}\right)\overline\delta.
\notag
\end{eqnarray}
Making $\Lambda_v\geq 1$ sufficiently large, so that
\begin{equation}
\frac{C_{v,e}}{\Lambda_v^{\frac1{2(1+\omega)}}}+\frac{C_e}{\Lambda_v^\frac12}\leq r_0
\label{ineq:condition-on-Lv-Lambdav-to-achieve-target-tensor-approximation}
\end{equation}
we can achieve 
\begin{eqnarray}
\left|\fint_{\mathbb T^d}\left(v_1\otimes v_1-v\otimes v-R_\ell\right)\,dx\right|
\leq r_0\overline\delta.
\label{ineq:target-tensor-approximation}
\end{eqnarray}

\nextpar{$C^0$-estimate on $\mathring R_1$}
We have
\begin{eqnarray}
\|\mathring R_1\|_0
\leq C_v\left(D\ell+\sqrt\delta D^\frac12\lambda^{2\alpha-\frac12}+\sqrt\delta D^\frac12\lambda^{\alpha-\frac12}\right)
\notag
\end{eqnarray}
using the fact that $\lambda\geq 1$ and thus we should keep the least negative of 
$\alpha-\frac12<2\alpha-\frac12<0$.
Note also that we have used that $\overline\delta\leq \sqrt\delta$.
Now (\ref{ineq:target-bound-on-|Ro1|0}) obtains provided
\begin{equation}
\frac{C_v}{L_v}+\frac{C_v}{\Lambda_v^\frac1{2(1+\varepsilon)}}\leq \eta.
\label{ineq:condition-on-Lv-Lambdav-to-achieve-target-|Ro1|0}
\end{equation}

\nextpar{$C^0$-estimate on $v_1-v$}
From Proposition~\ref{proposition:estimates-on-wo-wc-w},
\begin{equation}\|v_1-v\|_0=\|w\|_0\leq C_e\sqrt \delta\leq \frac M2\sqrt\delta
\label{ineq:condition-on-M-to-achieve-bound-on-|v1-v|0}
\end{equation}
by making $M$ sufficiently large. This is (\ref{ineq:target-bound-on-|v1-v|0}).
\bigskip

\nextpar{$C^0$-estimate on $p_1-p$}
The pressure $p_1$ has been defined in (\ref{eq:p1}) as $p_1=p+\tilde q-2\frac{(v-v_\ell)\cdot w}d$
where $\tilde q$ is given in (\ref{eq:q(tilde)}).
Making $M$ larger than previously if necessary (depending on $\nu$ in the case $d=2$),
we have 
$\|p_1-p\|_0\leq \frac{M^2}4\delta +\|v-v_\ell\|_0\|w\|_0$.
But from $\|v-v_\ell\|_0\leq CD\ell\leq C\overline\delta$,
we get $CD\ell C_e\sqrt \delta\leq CC_e\overline\delta\sqrt\delta\leq \frac{M^2}4\delta$.
Increasing $M$ if necessary, we get (\ref{ineq:target-bound-on-|p1-p|0}):
\begin{equation}
\|p_1-p\|_0\leq \frac{M^2}4\delta + C_e\delta\leq \frac{M^2}2\delta.
\label{ineq:condition-on-M-to-achieve-bound-on-|p1-p|0}
\end{equation}

\nextpar{$C^1$-estimates}
Since $\lambda\geq 1$ and $\alpha-\frac12<2\alpha-\frac12<0$, 
we have from Proposition~\ref{proposition:estimates-on-Ro1}
\begin{eqnarray}
\|\mathring R_1\|_{C^1}
\leq C_v\lambda\left\{\sqrt\delta\frac{\overline\delta}{L_v}
+\sqrt\delta D^\frac12\lambda^{2\alpha-\frac12}+\sqrt\delta D^\frac12\lambda^{\alpha-\frac12}\right\}
\leq
\lambda\overline\delta
\left(\frac{C_v}{L_v}+\frac{C_v}{\Lambda^\frac{1-\omega}{2(1+\omega)}}\right)
\notag
\end{eqnarray}
and therefore $\|\mathring R_1\|_1\leq\lambda\overline\delta$ provided
\begin{equation}
\frac{C_v}{L_v}+\frac{C_v}{\Lambda_v^\frac{1-\omega}{2(1+\omega)}}<1.
\label{ineq:condition-on-Lv-Lambdav-to-achieve-target-|Ro1|1}
\end{equation}

From Proposition~\ref{proposition:estimates-on-wo-wc-w},
$\|v_1\|_{C^1}\leq \|v\|_{C^1}+\|w\|_{C^1}\leq D+C_{e,v}\sqrt\delta\lambda$
so that, with $\overline\delta\leq \sqrt \delta$,
\begin{eqnarray}
\max\left\{\|v_1\|_{C^1},\|\mathring R_1\|_{C^1}\right\}
\leq 
 D+C_{e,v}\sqrt\delta\Lambda_v\left(\frac{D\delta}{\overline\delta^2}\right)^{1+\varepsilon}
\leq 2C_{e,v}\Lambda_v\delta^\frac32\left(\frac D{\overline\delta^2}\right)^{1+\varepsilon}
\notag
\end{eqnarray}
since $D\geq 1$ and $\delta^\frac32\geq \overline\delta^2$.
Now set $A:=2C_{e,v}\Lambda_v$.
From (\ref{ineq:a-priori-bound-on-|v(n)-v(0)|}),
we conclude
\begin{equation}A:=2C_{e}\Lambda_v.
\label{eq:A}
\end{equation}

\nextpar{Estimate on $\|v_1(\cdot,t)-v(\cdot,t)\|_{H^{-1}(\mathbb T^d)}$}
By construction we have $v_1-v=w=w_o+w_c$
and we will estimate $\|w_o\|_{H^{-1}(\mathbb T^d)}$ and $\|w_c\|_{H^{-1}(\mathbb T^d)}$ separately.

Let $f$ be any test vector field.
By definition (\ref{eq:wo(x,t)}) of $w_o$, 
and according to estimates from Propositions~\ref{proposition:stationary-phase-lemma}
and \ref{proposition:estimates-on-the-coefficients} in the Appendix, 
we have
\begin{eqnarray}
\left|\int_{\mathbb T^d}w_o\cdot f\,dx\right|
&\leq&C\sum_{|k|^2=\nu}\frac{\|a_k\|_1\|f\|_1}\lambda\notag\\
&\leq& \frac{C_e}{\Lambda_v^\frac12}D^{-\omega}\delta^{-\frac\omega{1-\omega}}\overline\delta^\frac{1+\omega}{1-\omega}\|f\|_1\notag\\
&\leq& \frac{C_e}{\Lambda_v^\frac12}\left(\frac{\overline\delta}\delta\right)^{\frac\omega{1-\omega}}\overline\delta^\frac1{1-\omega}\|f\|_1\notag\\
&\leq& \frac{C_e}{\Lambda_v^\frac12}\overline\delta\|f\|_1
\label{ineq:int(wo.f)}
\end{eqnarray}

Next, observe from (\ref{eq:properties-of-b(k)-in-3D})=(\ref{eq:properties-of-b(k)-in-2D}) and (\ref{eq:wo(x,t)})
that
\begin{eqnarray}
w_o(x,t)
=\frac1\lambda {\rm curl}\,\left(\sum_ka_k(x,t;\lambda t)\psi_k(\lambda x)\right)
-\frac1\lambda\sum_k\psi_k(\lambda x)\cdot {\rm curl}\,a_k(x,t;\lambda t)
\notag
\end{eqnarray}
and thus
$w_c(x,t)=-\mathcal Qw_o(x,t)
=\frac1\lambda \mathcal Q\left(\sum_{|k|^2=\nu} \psi_k(\lambda x)\cdot{\rm curl}\, a_k(x,t,\lambda t)\right)
=\frac1\lambda\mathcal Q\,u_c$.
(Recall the interpretation of the curl operator in dimension $d=2$ from Section~\ref{section:linear-spaces-of-stationary-flows}.)
The function $u_c$ is of the form
$u_c(x,t)=\sum_{|k|^2=\nu}\tilde c_k(x,t;\lambda t)e^{i\lambda k\cdot x}$
where the coefficients $\tilde c_k$ satisfy the same estimates as the coefficients $\nabla a_k$,
see Proposition~\ref{proposition:estimates-on-the-coefficients} in the Appendix.
Then, with $0<\gamma<1$ to be specified later,
we find
\begin{eqnarray}
\left|\int_{\mathbb T^d}w_c\cdot f\,dx\right|
&\leq& \frac1\lambda\left|\int_{\mathbb T^d}u_c\cdot \mathcal Qf\,dx\right|\notag\\
&\leq& C\frac1\lambda \sum_{|k|^2=\nu}\frac{\|\tilde c_k\|_\gamma\|\mathcal Q f\|_\gamma}{\lambda^\gamma}\notag\\
&\leq& \|\mathcal Qf\|_\gamma \frac1{\lambda^{1+\gamma}}C_e\sqrt\delta(\mu^{1+\gamma}D^{1+\gamma}+\mu D\ell^{-\gamma})\notag\\
&\leq&\|f\|_1C_e
\left\{
\lambda^{-\frac12(1+\gamma)}\sqrt\delta D^{\frac12(1+\gamma)}
+\lambda^{-\frac12-\gamma}\sqrt\delta D^{\frac12+\gamma}\overline\delta^{-\gamma}
\right\}\notag\\
&\leq&\|f\|_1C_e
\Bigg\{\Lambda_v^{-\frac12(1+\gamma)} D^{\frac12(1+\gamma)(1-\frac{1+\omega}{1-\omega})}\delta^{\frac12-\frac12(1+\gamma)\frac{1+\omega}{1-\omega}}\overline\delta^{(1+\gamma)\frac{1+\omega}{1-\omega}}\notag\\
&&\qquad +\Lambda_v^{-\frac12-\gamma} D^{(\frac12+\gamma)(1-\frac{1+\omega}{1-\omega})}\delta^{-(\frac12+\gamma)\frac{1+\omega}{1-\omega}+\frac12}\overline\delta^{(\frac12+\gamma)\frac{1+\omega}{1-\omega}-\gamma}
\Bigg\}\notag\\
&\leq&\|f\|_1C_e
\left\{\Lambda_v^{-\frac12(1+\gamma)}\overline\delta + \lambda^{-\frac12-\gamma}\overline\delta^{\frac12-\gamma}\right\}
\notag
\end{eqnarray}
where we have used that $\frac12(1+\gamma)\frac{1+\omega}{1-\omega}-\frac12>0$ for $0<\omega<1$ 
and $0<\gamma<1$, $D\geq 1$, and $\overline\delta\leq \delta$.
Fix $\gamma=\frac14$ so that
$\left|\int_{\mathbb T^d}w_c\cdot f\,dx\right|
\leq \|f\|_1\frac{C_e}{\Lambda_v^\frac58}\overline\delta^\frac14$.
From this and (\ref{ineq:int(wo.f)}) we conclude
$\left|\int_{\mathbb T^d}w\cdot f\,dx\right|\leq \frac{C_e}{\Lambda_v^\frac12}\overline\delta^\frac14\|f\|_1$
which implies that 
$\|w\|_{H^{-1}(\mathbb T^d)}\leq \frac{C_e}{\Lambda_v^\frac12}\overline\delta^\frac14$.
The bound (\ref{ineq:target-bound-on-|v1-v|H(-1)}) is satisfied 
provided
\begin{equation}
\frac{C_e}{\Lambda_v^\frac12}<r_0.
\label{ineq:condition-on-Lambdav-to-achieve-target-H(-1)-bound}
\end{equation}

\section{Proof of Theorem~\ref{theorem:weak-h-principle}, part 2}
\setcounter{counterpar}{0}
We consider again the sequence $(v_n,\mathring R_n,p_n)$ and the limit $v$
from Section~\ref{section:proof-of-part-1-of-theorems}.
We denote with some abuse $R_{n,\ell}=\rho_{n,\ell}{\rm Id}-\mathring R_{n,\ell}$ and $v_{n,\ell}$ the corresponding quantities
(since actually $\ell=\ell_n$).
Write
\begin{eqnarray}
v\otimes v-v_0\otimes v_0+\mathring R_0-\frac1d\Delta_0{\rm Id}
&=&\sum_{n=0}^\infty \left(v_{n+1}\otimes v_{n+1}-v_n\otimes v_n-R_{n,\ell}\right)\notag\\
&&+\sum_{n=0}^\infty \left(\rho_{n,\ell}-\frac1d\Delta_0(\delta_n-\delta_{n+1})\right){\rm Id}\notag\\
&&-\sum_{n=1}^\infty \mathring R_{n,\ell} + \mathring R_0-\mathring R_{0,\ell}.\notag
\end{eqnarray}
From (\ref{ineq:target-tensor-approximation}),
$$\left|\fint_{\mathbb T^d}\left(v_{n+1}\otimes v_{n+1}-v_n\otimes v_n-R_{n,\ell}\right)\,dx\right|
\leq r_0\delta_{n+1}\qquad (n\geq 0)$$
and from (\ref{ineq:n-th-step-bound-on-Ro}) and (\ref{ineq:|Ro(l)|<|Ro|}),
$$\left|\fint_{\mathbb T^d}\mathring R_{n,\ell}\,dx\right|\leq \eta\delta_{n-1}\qquad(n\geq 1).$$
As for the remaining term,
\begin{eqnarray}
d(2\pi)^d\rho_{n,\ell}
&=&e_0+\Delta_0(1-\delta_{n+1})-\fint_{\mathbb T^d}|v_{n,\ell}|^2\,dx\notag\\
&=&\Delta_0(\delta_n-\delta_{n+1})+\fint_{\mathbb T^d}\left(|v_n|^2-|v_{n,\ell}|^2\right)\,dx\notag\\
&&+e_0+\Delta_0(1-\delta_n)-\fint_{\mathbb T^d}|v_n|^2\,dx\notag
\end{eqnarray}
But from (\ref{ineq:n-th-step-energy-gap}) 
and from (\ref{ineq:int||v(l)|2-|v|2|})
we find
$$\left|d(2\pi)^d\rho_{n,\ell}-\Delta_0(\delta_n-\delta_{n+1})\right|
\leq \zeta/2\Delta_0\delta_n+C\eta\delta_n(\max e^\frac12+1).$$

Also, by definition (\ref{eq:ell}), $(D\ell)_0=\frac{\delta_1}{L_v}\leq\zeta$ and 
thus from (\ref{ineq:estimates-on-|v-v(l)|-and-|Ro-Ro(l)|}) we find
$$\|\mathring R_0-\mathring R_{0,\ell}\|\leq C\zeta.$$

With all the above we conclude
\begin{eqnarray}
&&\left|\fint_{\mathbb T^d}\left(v\otimes v-v_0\otimes v_0+\mathring R_0-\frac1d\Delta_0{\rm Id}\right)\,dx\right|\notag\\
&\leq&2\left(r_0+\zeta/2\Delta_0+C\frac{\Delta_0}dr_0(\max e^\frac12+1)\right)+C\zeta.
\label{ineq:condition-on-ro-zeta-eta-to-achieve-target-tensor-approximation}
\end{eqnarray}
Thus, making $r_0$ and $\zeta$ sufficiently small depending $\sigma, e$ and $\Delta_0$, 
we can achieve
$$\left|\fint_{\mathbb T^d} \left(v\otimes v-v_0\otimes v_0+\mathring R_0-\frac1d\Delta_0{\rm Id}\right)\,dx\right|<\sigma.$$
\stopproof

\section{Fixing the parameters $M$, $\zeta$,  $r_0$, $L_v$, and $\Lambda_v$}\label{section:fixing-the-parameters}
We list the requirements on the parameters $M$, $\zeta$, $r_0$, $L_v$, and $\Lambda_v$: 
\begin{itemize}
\item $M$: (\ref{ineq:condition-on-M-to-achieve-bound-on-|v1-v|0}),
(\ref{ineq:condition-on-M-to-achieve-bound-on-|p1-p|0});
\item $\zeta$: (\ref{ineq:condition-on-zeta-so-that-initial-R-can-be-generated}),
(\ref{ineq:condition-on-ro-zeta-eta-to-achieve-target-tensor-approximation});
\item $r_0$: (\ref{ineq:condition-on-r0-for-H(-1)-approximation}),
(\ref{ineq:int||v(l)|2-|v|2|}),
(\ref{ineq:condition-on-eta-for-lower-bound-on-rho(l)}), 
(\ref{ineq:condition-on-eta-so-that-R(l)/rho(l)-is-in-N}),
(\ref{ineq:condition-on-ro-zeta-eta-to-achieve-target-tensor-approximation});
\item $L_v$: (\ref{ineq:condition-on-Lv-to-achieve-conditions-on-parameters}),
(\ref{ineq:condition-on-Lv-Lambdav-to-achieve-conditions-on-parameters}),
(\ref{ineq:condition-on-Lv-Lambdav-to-achieve-target-energy-gap}), 
(\ref{ineq:condition-on-Lv-Lambdav-to-achieve-target-tensor-approximation}),
(\ref{ineq:condition-on-Lv-Lambdav-to-achieve-target-|Ro1|0}),
(\ref{ineq:condition-on-Lv-Lambdav-to-achieve-target-|Ro1|1});
\item $\Lambda_v$: (\ref{ineq:condition-on-Lv-Lambdav-to-achieve-conditions-on-parameters}),
(\ref{ineq:condition-on-Lv-Lambdav-to-achieve-target-energy-gap}),
(\ref{ineq:condition-on-Lv-Lambdav-to-achieve-target-tensor-approximation}),
(\ref{ineq:condition-on-Lv-Lambdav-to-achieve-target-|Ro1|0})
(\ref{ineq:condition-on-Lv-Lambdav-to-achieve-target-|Ro1|1}),
(\ref{ineq:condition-on-Lambdav-to-achieve-target-H(-1)-bound}).
\end{itemize}
Recall that $d, e, \Delta_0, (v_0,\mathring R_0,p_0)$, and $\varepsilon$ (hence $\omega$) are given.
We set the parameters $M, \zeta, r_0, L_v$ and $\Lambda_v$ in this order as follows:
\begin{enumerate}
\item set $M$ larger than a constant $C_e$ so that it satisfies
(\ref{ineq:condition-on-M-to-achieve-bound-on-|v1-v|0}),
(\ref{ineq:condition-on-M-to-achieve-bound-on-|p1-p|0})
(specifically, the constant depends on $e$ and $\sup_{k,y,\tau}\|a_k(\cdot,y;\tau)\|_1$);
\item set $\zeta$ sufficiently small depending on $(v_0,\mathring R_0,p_0)$ and $\sigma$
so that it satisfies
(\ref{ineq:condition-on-zeta-so-that-initial-R-can-be-generated}),
(\ref{ineq:condition-on-ro-zeta-eta-to-achieve-target-tensor-approximation});
\item set $r_0$ sufficiently small depending on $d,\sigma, \Delta_0$ and $e$
so that it satisfies
(\ref{ineq:condition-on-r0-for-H(-1)-approximation}),
(\ref{ineq:int||v(l)|2-|v|2|}),
(\ref{ineq:condition-on-eta-for-lower-bound-on-rho(l)}), 
(\ref{ineq:condition-on-eta-so-that-R(l)/rho(l)-is-in-N}),
(\ref{ineq:condition-on-ro-zeta-eta-to-achieve-target-tensor-approximation});
\item set $L_v\geq 1$ sufficiently large 
depending on $r_0$, $\Delta_0$, and $\varepsilon$ 
so that it satisfies
(\ref{ineq:condition-on-Lv-to-achieve-conditions-on-parameters}),
(\ref{ineq:condition-on-Lv-Lambdav-to-achieve-target-energy-gap}), 
(\ref{ineq:condition-on-Lv-Lambdav-to-achieve-target-tensor-approximation}),
(\ref{ineq:condition-on-Lv-Lambdav-to-achieve-target-|Ro1|0}),
(\ref{ineq:condition-on-Lv-Lambdav-to-achieve-target-|Ro1|1});
\item finally set $\Lambda_v$ sufficiently large, depending on $\zeta, \Delta_0, \varepsilon,r_0$, and $L_v$, 
so that it satisfies
(\ref{ineq:condition-on-Lv-Lambdav-to-achieve-conditions-on-parameters}),
(\ref{ineq:condition-on-Lv-Lambdav-to-achieve-target-energy-gap}),
(\ref{ineq:condition-on-Lv-Lambdav-to-achieve-target-tensor-approximation}),
(\ref{ineq:condition-on-Lv-Lambdav-to-achieve-target-|Ro1|0})
(\ref{ineq:condition-on-Lv-Lambdav-to-achieve-target-|Ro1|1}),
(\ref{ineq:condition-on-Lambdav-to-achieve-target-H(-1)-bound}).
\end{enumerate}

\appendix

\section{Estimates from \cite{DS-Hoelder}}
\setcounter{counterpar}{0}
The following Propositions have been proved in \cite{DS-Hoelder} and \cite{DS3}.

\begin{proposition}[Schauder estimates for elliptic operators]\label{proposition:Schauder-estimates}
Let $d\geq 2$.
For any $\alpha\in (0,1)$ and any $m\in\mathbb N$ there exists a constant $C_s(d,m,\alpha)$
so that the following estimates hold.
\begin{eqnarray}
\|\mathcal Qv\|_{m+\alpha}&\leq&C_s(d,m,\alpha)\|v\|_{m+\alpha}\notag\\
\|\mathcal Pv\|_{m+\alpha}&\leq&C_s(d,m,\alpha)\|v\|_{m+\alpha}\notag\\
\|\mathcal Rv\|_{m+1+\alpha}&\leq&C_s(d,m,\alpha)\|v\|_{m+\alpha}\notag\\
\|\mathcal R\,{\rm div}\,A\|_{m+\alpha}&\leq&C_s(d,m,\alpha)\|A\|_{m+\alpha}\notag\\
\|\mathcal R\mathcal Q{\rm div}\,A\|_{m+\alpha}&\leq&C_s(d,m,\alpha)\|A\|_{m+\alpha}\notag
\end{eqnarray}
\end{proposition}
\startproof
This is Proposition~4.3 of \cite{DS-Hoelder}, valid as it is for any $d\geq 2$
provided $\mathcal R$ is defined according to Definition~\ref{definition:the-operator-R}.
\stopproof

\begin{proposition}[Stationary phase lemma]\label{proposition:stationary-phase-lemma}
Let $d\geq 1$.
For $k\in\mathbb Z^d$ and $\lambda\geq 1$,
\begin{enumerate}
\item For any function $a\in C^\infty(\mathbb T^d)$ and $m\in\mathbb N$ we have
$$\left|\int_{\mathbb T^d}a(x)e^{i\lambda k\cdot x}\,dx\right|\leq\frac{[a]_m}{\lambda^m}.$$
\item Let $k\in\mathbb Z^d\setminus \{0\}$.
For any vector field $F\in C^\infty(\mathbb T^d,\mathbb R^d)$ let $F_\lambda(x):=F(x)e^{i\lambda k\cdot x}$.
Then,
\begin{eqnarray}
\|\mathcal R\,F_\lambda\|_\alpha&\leq&
\frac{C_s}{\lambda^{1-\alpha}}\|F\|_0+\frac{C_s}{\lambda^{m-\alpha}}\left[F\right]_m\frac{C_s}{\lambda^m}\left[F\right]_{m+\alpha}\notag\\
\|\mathcal R\,\mathcal Q\,F_\lambda\|_\alpha&\leq&
\frac{C_s}{\lambda^{1-\alpha}}\|F\|_0+\frac{C_s}{\lambda^{m-\alpha}}\left[F\right]_m\frac{C_s}{\lambda^m}\left[F\right]_{m+\alpha}\notag
\end{eqnarray}
where $C_s=C_s(d,m,\alpha)$ ({\it i.e.} they do not depend on $\lambda$ nor $k\neq 0$).
\end{enumerate}
\end{proposition}
\startproof This is Proposition~4.4 of \cite{DS-Hoelder}, valid as it is for any $d\geq 1$.
\stopproof

Of the estimates from Proposition~5.1 from \cite{DS-Hoelder}, we will only recall the one
which is explicitly used here (in the estimate of $\|w_c\|_{H^{-1}(\mathbb T^d)}$).
\begin{proposition}[Estimates on the coefficients]\label{proposition:estimates-on-the-coefficients}
Let $a_k\in C^\infty(\mathbb T^d\times\mathbb S^1\times \mathbb R)$ be given by (\ref{eq:ak}).
For any $r\geq 1$,
$$\|a_k(\cdot, s;\tau)\|_r\leq C_e\sqrt\delta\left(\mu^rD^r+\mu D\ell^{1-r}\right).
$$
\end{proposition}

\begin{proposition}[Estimates on $w_o$, $w_c$, and $w$]\label{proposition:estimates-on-wo-wc-w}
Assuming (\ref{ineq:conditions-on-the-parameters}) and $r\geq 0$, we have
$$\|w_o\|_r\leq C_e\sqrt \delta\lambda^r$$
$$\|\partial_tw_o\|_r\leq C_v\sqrt \delta \lambda^{r+1}$$
and for $r>0$, $r\not\in\mathbb N$,
$$\|w_c\|_r\leq C_{e,s}\sqrt\delta D\mu\lambda^{r-1}$$
$$\|\partial_tw_c\|_r\leq C_{v,s}\sqrt \delta D\mu\lambda^r.$$
In particular
$$\|w\|_0\leq C_e\sqrt \delta$$
$$\|w\|_1\leq C_{e,v}\sqrt \delta \lambda.$$
\end{proposition}
\startproof
This is Proposition~6.1 of \cite{DS-Hoelder}.
\stopproof

\begin{proposition}[Estimates on the energy]\label{proposition:estimates-on-the-energy}
For any $\alpha\in(0,\frac\omega{1+\omega})$ 
there is a constant $C_{v,s}$ depending only on $d$, $\alpha$, $e$, and $\|v\|_0$,
such that,
under the assumptions (\ref{ineq:conditions-on-the-parameters}), we have
$$\left|e_0(t)+\Delta_0(t)(1-\overline\delta)-\fint_{\mathbb T^d}|v_1(x,t)|^2\,dx\right|
\leq C_eD\ell+C_{v,s}\sqrt \delta\mu D\lambda^{\alpha-1}$$
and 
$$\left|\fint_{\mathbb T^d}\left(v_1\otimes v_1-v\otimes v-R_\ell\right)\,dx\right|
\leq C_{v,e} \sqrt\delta D\mu\lambda^{\alpha-1}+C_e\delta D\mu\lambda^{-1}.
$$
\end{proposition}
\startproof
The first estimate is Proposition~7.1 of \cite{DS-Hoelder}.
The second estimate holds since 
since the first and second term of
$$\int_{\mathbb T^d}\left(v_1\otimes v_1-v\otimes v-R_\ell\right)\,dx
=\int_{\mathbb T^d}\left(v_1\otimes v_1-v\otimes v-w_o\otimes w_o\right)\,dx
+\int_{\mathbb T^d}\left(w_o\otimes w_o-R_\ell\right)\,dx.
$$
are estimated exactly as $\int_{\mathbb T^d}\left(|v_1|^2-|v|^2-|w_o|^2\right)\,dx$
and $\int_{\mathbb T^d}\left(|w_o|^2-{\rm tr}\,R_\ell\right)\,dx$
in the proof of Proposition~7.1 of \cite{DS-Hoelder}, respectively.
\stopproof

\begin{proposition}[Estimates on $\mathring R_1$]\label{proposition:estimates-on-Ro1}
For every $\alpha\in(0,\frac\omega{1+\omega})$ there is a constant $C_{v,s}$ depending only
on $d, \alpha, \omega, e$, and $\|v\|_0$, such that,
under the assumptions (\ref{ineq:conditions-on-the-parameters}), we have
$$\|\mathring R_1\|_0\leq C_{v,s}\left(D\ell+\sqrt \delta D\mu\lambda^{2\alpha-1}+\sqrt\delta\mu^{-1}\lambda^\alpha\right)$$
$$\|\mathring R_1\|_1
\leq C_{v,s}\lambda\left(\sqrt\delta D\ell+\sqrt\delta D\mu\lambda^{2\alpha-1}+\sqrt\delta\mu^{-1}\lambda^\alpha\right)$$
\end{proposition}
\startproof
This is Proposition~8.1 of \cite{DS-Hoelder}.
For convenience for the reader, we recall briefly how the ``zero-mode'' of $w_o\otimes w_o$ cancels
with $R$.
From definition (\ref{eq:W(y,s,tau,xi)}) of $W$ we have
$$W\otimes W(y,s;\tau,\xi)= U_o(y,s)+\sum_{1\leq |k|\leq 2\nu} U_k(y,s;\tau)e^{ik\cdot\xi}
$$
for some coefficients $U_k$.
The ``zero-mode'' $U_0(y,s)$ is precisely $R_\ell(y,s)$
since
\begin{eqnarray}
\fint_{\mathbb T^d}W\otimes W\,d\xi
&\stackrel{}=&\rho_\ell\sum_{j=1}^{2^d}\sum_{k\in\Lambda_j}\left(\gamma_k\left(\frac{R_\ell}{\rho_\ell}\right)\right)^2|\phi_k^{(j)}(v,\tau)|^2\left({\rm Id}-\frac k{|k|}\otimes \frac k{|k|}\right)\notag\\
&\stackrel{}=&\rho_\ell\sum_{j=1}^{2^d}\sum_{k\in\Lambda_j}\sum_{l\in\mathcal C_j}\left(\gamma_k\left(\frac{R_\ell}{\rho_\ell}\right)\right)^2\alpha_l^2(v)\left({\rm Id}-\frac k{|k|}\otimes \frac k{|k|}\right)\notag\\
&\stackrel{}=&R_\ell\sum_{j=1}^{2^d}\sum_{l\in\mathcal C_j}\alpha_l^2(v)\notag\\
&\stackrel{}=&R_\ell.\notag
\end{eqnarray}
The crucial identities are (\ref{eq:ave(WxW)=sum|ak|2Mk}),
(\ref{eq:R=sum-of-Mk}),
(\ref{eq:sum-of-ak-squared=1}), and (\ref{eq:|phi(j)k|2=sum-of-alpha(k)2}).
Thus,
\begin{eqnarray}
{\rm div}\,(w_o\otimes w_0+\mathring R_\ell+\tilde q{\rm Id})
&=&{\rm div}_y\,(W\otimes W-R_\ell)+\nabla_yQ +\lambda {\rm div}_\xi\,(W\otimes W+Q{\rm Id})\notag\\
&=&\sum_{1\leq |k|\leq 2\nu}\left({\rm div}_y\,U_k+\nabla_y\tilde a_k{\rm Id}\right)e^{i\lambda k\cdot x}
\notag
\end{eqnarray}
where a cancelation occurs since $(W,Q)$ is a stationary solution (in the $\xi$-variable) to the Euler equations.
Here we have used that $\rho_\ell=\rho_\ell(t)$,
and the $\tilde a_k$'s are the coefficients of $Q$, see (\ref{eq:Q}).
In the end, ${\rm div}\,(w_o\otimes w_o+R_\ell+\tilde q{\rm Id})$ is oscillatory.

The other terms in ${\rm div}\,\mathring R_1$ are linear in $w$ and hence are also oscillatory.
\stopproof


\bibliographystyle{amsalpha}
\bibliography{eulerbib}

\providecommand{\bysame}{\leavevmode\hbox to3em{\hrulefill}\thinspace}
\providecommand{\MR}{\relax\ifhmode\unskip\space\fi MR }
\providecommand{\MRhref}[2]{%
  \href{http://www.ams.org/mathscinet-getitem?mr=#1}{#2}
}
\providecommand{\href}[2]{#2}
\begin{thebibliography}{CDLS11}

\bibitem[Bor65]{Borisov65}
Ju.~F. Borisov, \emph{{$C^{1,\alpha}$-isometric immersions of {Riemannian}
  spaces}}, Doklady \textbf{163} (1965), 869--871.

\bibitem[Bor04]{Borisov2004}
Yu.F. Borisov, \emph{{Irregular $C^{1,\beta}$-surfaces with analytic metric.}},
  Sib. Mat. Zh. \textbf{45} (2004), no.~1, 25--61 (Russian, English).

\bibitem[CDLS11]{CDSz}
Sergio Conti, Camillo De~Lellis, and L{\'a}szl{\'o} Sz{\'e}kelyhidi, Jr.,
  \emph{$h$-principle and rigidity for {$C^{1,\alpha}$} isometric embeddings},
  To appear in the Proceedings of the Abel Symposium 2010 (2011).

\bibitem[CET94]{ConstantinETiti}
Peter Constantin, Weinan E, and Edriss~S. Titi, \emph{Onsager's conjecture on
  the energy conservation for solutions of {E}uler's equation}, Comm. Math.
  Phys. \textbf{165} (1994), no.~1, 207--209. \MR{MR1298949 (96e:76025)}

\bibitem[CFG11]{CFG}
Diego Cordoba, Daniel Faraco, and Francisco Gancedo, \emph{Lack of uniqueness
  for weak solutions of the incompressible porous media equation}, Arch.
  Ration. Mech. Anal. \textbf{200} (2011), no.~3, 725--746. \MR{2796131}

\bibitem[Chi12]{Chiodaroli}
E.~Chiodaroli, \emph{A counterexample to well-posedness of entropy solutions to
  the compressible {E}uler system}, Preprint (2012).

\bibitem[DLS09]{DS1}
Camillo De~Lellis and L{\'a}szl{\'o} Sz{\'e}kelyhidi, Jr., \emph{The {E}uler
  equations as a differential inclusion}, Ann. of Math. (2) \textbf{170}
  (2009), no.~3, 1417--1436. \MR{2600877 (2011e:35287)}

\bibitem[DLS10]{DS2}
\bysame, \emph{On admissibility criteria for weak solutions of the {E}uler
  equations}, Arch. Ration. Mech. Anal. \textbf{195} (2010), no.~1, 225--260.
  \MR{2564474 (2011d:35386)}

\bibitem[DLS12a]{DS3}
\bysame, \emph{Continuous dissipative {E}uler flows}, preprint, {\texttt
  arxiv.1202.1751}, (2012).

\bibitem[DLS12b]{DS-Hoelder}
\bysame, \emph{Dissipative {E}uler flows and {O}nsager's conjecture}, preprint,
  {\texttt arxiv.1205.3626}, (2012).

\bibitem[DLS12c]{DSsurvey}
\bysame, \emph{The $h$-principle and the equations of fluid dynamics}, Bull.
  Amer. Math. Soc. (N.S.) \textbf{49} (2012), no.~3, 347--375.

\bibitem[EM02]{Eliashberg}
Y.~Eliashberg and N.~Mishachev, \emph{Introduction to the {$h$}-principle},
  Graduate Studies in Mathematics, vol.~48, American Mathematical Society,
  Providence, RI, 2002. \MR{MR1909245 (2003g:53164)}

\bibitem[Eyi94]{Eyink}
Gregory~L. Eyink, \emph{Energy dissipation without viscosity in ideal
  hydrodynamics. {I}. {F}ourier analysis and local energy transfer}, Phys. D
  \textbf{78} (1994), no.~3-4, 222--240. \MR{MR1302409 (95m:76020)}

\bibitem[Fri95]{FrischBook}
Uriel Frisch, \emph{Turbulence}, Cambridge University Press, Cambridge, 1995,
  The legacy of A. N. Kolmogorov. \MR{MR1428905 (98e:76002)}

\bibitem[Gro86]{Gromov}
Mikhael Gromov, \emph{Partial differential relations}, Ergebnisse der
  Mathematik und ihrer Grenzgebiete (3), vol.~9, Springer-Verlag, Berlin, 1986.
  \MR{90a:58201}

\bibitem[Ham82]{Hamilton:Nash-Moser-IFT}
R.~S. Hamilton, \emph{The {I}nverse {F}unction {T}heorem of {N}ash and
  {M}oser}, Bull. Amer. Math. Soc. (N.S.) \textbf{7} (1982), no.~1, 65--222.

\bibitem[Maj91]{MAJDA:nonlinear-analysis-applied-math}
Andrew Majda, \emph{The interaction of nonlinear analysis and modern applied
  mathematics}, 175--191.

\bibitem[MB02]{MajdaBook}
Andrew~J. Majda and Andrea~L. Bertozzi, \emph{Vorticity and incompressible
  flow}, Cambridge Texts in Applied Mathematics, vol.~27, Cambridge University
  Press, Cambridge, 2002. \MR{MR1867882 (2003a:76002)}

\bibitem[Nas54]{Nash54}
J.~Nash, \emph{{$C^1$} isometric imbeddings}, Ann. Math. \textbf{60} (1954),
  383--396.

\bibitem[Ons49]{Onsager}
L.~Onsager, \emph{Statistical hydrodynamics}, Nuovo Cimento (9) \textbf{6}
  (1949), no.~Supplemento, 2(Convegno Internazionale di Meccanica Statistica),
  279--287. \MR{MR0036116 (12,60f)}

\bibitem[Sch93]{Scheffer93}
Vladimir Scheffer, \emph{An inviscid flow with compact support in space-time},
  J. Geom. Anal. \textbf{3} (1993), no.~4, 343--401. \MR{MR1231007 (94h:35215)}

\bibitem[Shn97]{Shnirelman1}
A.~Shnirelman, \emph{On the nonuniqueness of weak solution of the {E}uler
  equation}, Comm. Pure Appl. Math. \textbf{50} (1997), no.~12, 1261--1286.
  \MR{MR1476315 (98j:35149)}

\bibitem[Shn00]{Shnirelmandecrease}
\bysame, \emph{Weak solutions with decreasing energy of incompressible {E}uler
  equations}, Comm. Math. Phys. \textbf{210} (2000), no.~3, 541--603.
  \MR{MR1777341 (2002g:76009)}

\bibitem[Shv11]{Shvydkoy}
R.~Shvydkoy, \emph{Convex integration for a class of active scalar equations},
  J. Amer. Math. Soc. \textbf{24} (2011), no.~4, 1159--1174. \MR{2813340}

\bibitem[Spr98]{SPRING:book}
D.~Spring, \emph{{Convex integration theory. Solutions to the h-principle in
  geometry and topology}}, {Birkh\"auser Verlag}, 1998.

\bibitem[SW11]{SzWie}
L{\'a}szl{\'o} Sz{\'e}kelyhidi, Jr. and E.~Wiedemann, \emph{Young measures
  generated by ideal incompressible fluid flows}, Preprint (2011).

\bibitem[Sz{\'e}11]{Szekelyhidi}
L{\'a}szl{\'o} Sz{\'e}kelyhidi, Jr., \emph{Relaxation of the incompressible
  porous medium equation}, Preprint (2011).

\bibitem[Wie11]{Wiedemann}
E.~Wiedemann, \emph{Existence of weak solutions for the incompressible {E}uler
  equations}, Ann. Inst. H. Poincar\'e Anal. Non Lin\'eaire \textbf{28} (2011),
  no.~5, 727--730.

\end{thebibliography}

\end{document}